\documentclass[12pt]{amsart}
\usepackage{amsmath,amsthm}
\usepackage{amssymb}
\usepackage[all]{xypic}
\SelectTips{cm}{12}
\UseTips{}
\usepackage{euscript}
\begin{document} 
%
%
\theoremstyle{plain}
\swapnumbers
	\newtheorem{thm}{Theorem}[section]
	\newtheorem{prop}[thm]{Proposition}
	\newtheorem{lemma}[thm]{Lemma}
	\newtheorem{cor}[thm]{Corollary}
	\newtheorem{fact}[thm]{Fact}
	\newtheorem{subsec}[thm]{}
\theoremstyle{definition}
	\newtheorem{assume}[thm]{Assumption}
	\newtheorem{defn}[thm]{Definition}
	\newtheorem{example}[thm]{Example}
	\newtheorem{examples}[thm]{Examples}
	\newtheorem{claim}[thm]{Claim}
	\newtheorem{notn}[thm]{Notation}
	\newtheorem{construct}[thm]{Construction}
\theoremstyle{remark}
        \newtheorem{remark}[thm]{Remark}
        \newtheorem{remarks}[thm]{Remarks}
	\newtheorem{ack}[thm]{Acknowledgements}
\newenvironment{myeq}[1][]
{\stepcounter{thm}\begin{equation}\tag{\thethm}{#1}}
{\end{equation}}
\newcommand{\mydiag}[2][]{\myeq[#1]{\xymatrix{#2}}}
\newcommand{\mydiagram}[2][]
{\stepcounter{thm}\begin{equation}
     \tag{\thethm}{#1}\vcenter{\xymatrix{#2}}\end{equation}}
%
\newenvironment{mysubsection}[2][]
{\begin{subsec}\begin{upshape}\begin{bfseries}{#2.}
\end{bfseries}{#1}}
{\end{upshape}\end{subsec}}
\newenvironment{mysubsect}[2][]
{\begin{subsec}\begin{upshape}\begin{bfseries}{#2\vsm.}
\end{bfseries}{#1}}
{\end{upshape}\end{subsec}}
\newcommand{\sect}{\setcounter{thm}{0}\section}
\newcommand{\bproof}{\noindent\textbf{Proof:}\ }
\newcommand{\wh}{\ -- \ }
\newcommand{\w}[2][ ]{\ \ensuremath{#2}{#1}\ }
\newcommand{\ww}[1]{\ \ensuremath{#1}}
\newcommand{\wb}[2][ ]{\ (\ensuremath{#2}){#1}\ }
\newcommand{\wref}[2][ ]{\ \eqref{#2}{#1}\ }
%
%
\newcommand{\xra}[1]{\xrightarrow{#1}}
\newcommand{\xla}[1]{\xleftarrow{#1}}
\newcommand{\Ra}{\Rightarrow}
\newcommand{\hra}{\hookrightarrow}
\newcommand{\bstar}{\mbox{\large $\star$}}
\newcommand{\adj}[2]{\substack{{#1}\\ \rightleftharpoons \\ {#2}}}
\newcommand{\hsp}{\hspace{10 mm}}
\newcommand{\hs}{\hspace{5 mm}}
\newcommand{\hsm}{\hspace{2 mm}}
\newcommand{\vs}{\vspace{7 mm}}
\newcommand{\vsm}{\vspace{2 mm}}
\newcommand{\rest}[1]{\lvert_{#1}}
\newcommand{\lra}[1]{\langle{#1}\rangle}
\newcommand{\lin}[1]{\{{#1}\}}
\newcommand{\EQUIV}{\Leftrightarrow}
\newcommand{\epic}{\to\hspace{-5 mm}\to}
\newcommand{\xepic}[1]{\xrightarrow{#1}\hspace{-5 mm}\to}
\newcommand{\hotimes}{\hat{\otimes}}
\newcommand{\hy}[2]{{#1}\text{-}{#2}}
%
%
\newcommand{\co}[1]{c({#1})}
\newcommand{\csk}[1]{\operatorname{csk}_{#1}}
\newcommand{\sk}[1]{\operatorname{sk}_{#1}}
\newcommand{\q}[1]{^{({#1})}}
\newcommand{\li}[1]{_{({#1})}}
\newcommand{\rz}{\rho_{0}}
%
%
\newcommand{\ab}{\operatorname{ab}}
\newcommand{\Arr}{\operatorname{Arr}}
\newcommand{\Aut}{\operatorname{Aut}}
\newcommand{\Card}{\operatorname{Card}}
\newcommand{\Cok}{\operatorname{Coker}\,}
\newcommand{\cof}{\operatorname{cof}}
\newcommand{\cf}{\operatorname{cf}}
\newcommand{\Coef}{\operatorname{Coef}}
\newcommand{\colim}{\operatorname{colim}}
\newcommand{\diag}{\operatorname{diag}}
\newcommand{\Ext}{\operatorname{Ext}}
\newcommand{\Fac}{\operatorname{Fac}}
\newcommand{\hc}[1]{\operatorname{hc}_{#1}}
\newcommand{\ho}{\operatorname{ho}}
\newcommand{\Hom}{\operatorname{Hom}}
\newcommand{\Id}{\operatorname{Id}}
\newcommand{\inc}{\operatorname{inc}}
\newcommand{\Ker}{\operatorname{Ker}\,}
\newcommand{\map}{\operatorname{map}}
\newcommand{\mapa}{\map_{\ast}}
\newcommand{\Obj}{\operatorname{Obj}}
\newcommand{\op}{\sp{\operatorname{op}}}
\newcommand{\red}{\operatorname{red}}
\newcommand{\Stov}{\operatorname{St}}
\newcommand{\Tor}{\operatorname{Tor}}
\newcommand{\we}{\operatorname{w.e.}}
\newcommand{\sss}{\hspace*{1 mm}\sp{s}}
%
%
\newcommand{\A}{{\EuScript A}}
\newcommand{\tA}{\tilde{\A}}
\newcommand{\hA}{\hat  {\A}}
\newcommand{\Aj}[1]{A\li{#1}}
\newcommand{\B}{{\mathcal B}}
\newcommand{\hB}{\hat{\B}}
\newcommand{\tB}{\tilde{\B}}
\newcommand{\C}{{\mathcal C}}
\newcommand{\tCA}{\tilde{\C}_{\A}}
\newcommand{\sC}{\sss\C}
\newcommand{\CA}{\C_{\A}}
\newcommand{\sCA}{\sp{s}\CA}
\newcommand{\CAst}{\CA^{\Stov}}
\newcommand{\CAs}{\sss\CAst}
\newcommand{\uCAs}{\CA^{\Stov}}
\newcommand{\CAm}{\sCA^{\min}}
\newcommand{\CAp}{\sCA\sp{+}}
\newcommand{\CAsp}{(\CAs)\sp{+}}
\newcommand{\CAstp}{(\CAst)\sp{+}}
\newcommand{\CL}{\CAp/\Lambda}
\newcommand{\hCA}{\sss\CA'}
\newcommand{\CG}{\C_{\Gamma}}
\newcommand{\CX}{\sss\C_{X}}
\newcommand{\CXp}{\sss\C_{X'}}
\newcommand{\D}{{\mathcal D}}
\newcommand{\tD}{\tilde{\D}}
\newcommand{\tDp}{\tilde{\D}\sp{+}}
\newcommand{\Dpp}{\D\sp{+}}
\newcommand{\hD}{\hat{\D}}
\newcommand{\E}{{\mathcal E}}
\newcommand{\Ep}{\E\sp{+}}
\newcommand{\F}{{\mathcal F}}
\newcommand{\hF}{\hat{F}}
\newcommand{\tF}{\tilde{F}}
\newcommand{\FA}{\F_{\A}}
\newcommand{\FAp}{\FA\sp{+}}
\newcommand{\Fs}{\F_{\bullet}}
\newcommand{\G}{{\mathcal G}}
\newcommand{\K}{{\mathcal K}}
\newcommand{\tK}{\tilde{K}}
\newcommand{\bK}{\bar{K}}
\newcommand{\LL}{{\mathcal L}}
\newcommand{\tL}{\tilde{\LL}}
\newcommand{\M}{{\mathcal M}}
\newcommand{\MA}{\M\sb{\A}}
\newcommand{\MAC}{\MA\sp{\CA'}}
\newcommand{\MAn}[1]{\M\sb{\An{#1}}}
\newcommand{\MAm}{\MA\sp{\min}}
\newcommand{\MAs}{\MA\sp{\Stov}}
\newcommand{\N}{{\mathcal N}}
\newcommand{\tN}{\tilde{N}}
\newcommand{\OO}{{\mathcal O}}
\newcommand{\Op}{\OO\sp{+}}
\newcommand{\cP}{{\mathcal P}}
\newcommand{\hQ}{\hat{Q}}
\newcommand{\R}[1]{{\mathcal R}_{#1}}
\newcommand{\Ss}{{\mathcal S}}
\newcommand{\Sf}{\Ss_{f}}
\newcommand{\Sa}{\Ss_{\ast}}
\newcommand{\Sr}{\Ss_{\red}}
\newcommand{\TT}{{\mathcal T}}
\newcommand{\Ta}{\TT_{\ast}}
\newcommand{\TA}{T_{A}}
\newcommand{\TAi}[2]{T_{#1}\q{#2}}
\newcommand{\UA}[1]{U_{#1}}
\newcommand{\U}{{\mathcal U}}
\newcommand{\V}{{\mathcal V}}
\newcommand{\W}{{\mathcal W}}
\newcommand{\X}{{\mathcal X}}
\newcommand{\Y}{{\mathcal Y}}
\newcommand{\Z}{{\mathcal Z}}
%
%
\newcommand{\Alg}[1]{{#1}\text{-}{\EuScript Alg}}
\newcommand{\Ab}{{\EuScript Ab}}
\newcommand{\Abgp}{{\Ab\Gp}}
\newcommand{\Cat}{{\EuScript Cat}}
\newcommand{\iO}[1]{({#1},\OO)}
\newcommand{\iOp}[1]{({#1},\Op)}
\newcommand{\CO}{\iO{\Set^{\Box}}}
\newcommand{\GO}{(\Gpd,\OO)}
\newcommand{\GOC}{\hy{\GO}{\Cat}}
\newcommand{\SO}{\iO{\Ss}}
\newcommand{\SOp}{(\Ss,\Op)}
\newcommand{\SaO}{\iO{\Sa}}
\newcommand{\SaOp}{\iOp{\Sa}}
\newcommand{\OC}{\hy{\OO}{\Cat}}
\newcommand{\OpC}{\hy{\Op}{\Cat}}
\newcommand{\iOC}[1]{\hy{\iO{#1}}{\Cat}}
\newcommand{\iOpC}[1]{\hy{\iOp{#1}}{\Cat}}
\newcommand{\SOC}{\iOC{\Ss}}
\newcommand{\SaOC}{\iOC{\Sa}}
\newcommand{\SaOpC}{\iOpC{\Sa}}
\newcommand{\COC}{\iOC{\Set^{\Box}}}
\newcommand{\VO}{\iO{\V}}
\newcommand{\VOC}{\iOC{\V}}
\newcommand{\VOpC}{\iOpC{\V}}
\newcommand{\Ch}{{\EuScript Chain}}
\newcommand{\DG}{{\EuScript D}i{\EuScript G}}
\newcommand{\Gp}{{\EuScript Gp}}
\newcommand{\Gpd}{{\EuScript Gpd}}
\newcommand{\NS}{{\EuScript NS}}
\newcommand{\Set}{{\EuScript Set}}
\newcommand{\Seta}{\Set_{\ast}}
\newcommand{\Track}{{\EuScript Trk}}
\newcommand{\RM}[1]{{#1}\text{-}{\EuScript Mod}}
%
%
\newcommand{\ma}{mapping algebra}
\newcommand{\Ama}{$\A$-\ma}
\newcommand{\Sma}{Stover mapping algebra}
\newcommand{\Pal}{$\Pi$-algebra}
\newcommand{\PiA}{\Pi_{\A}}
\newcommand{\PuA}{\Pi^{\A}}
\newcommand{\PiC}{\Pi_{\C}}
\newcommand{\PiD}{\Pi_{\D}}
\newcommand{\PAa}{$\PiA$-algebra}
\newcommand{\PAlg}{\Alg{\Pi}}
\newcommand{\TAlg}{\Alg{\Theta}}
\newcommand{\PAAlg}{\Alg{\PiA}}
\newcommand{\Tal}{$\Theta$-algebra}
%
%
\newcommand{\fG}{\mathfrak{G}}
%
%
\newcommand{\bF}{\mathbb F}
\newcommand{\Fp}{\bF_{p}}
\newcommand{\bN}{\mathbb N}
\newcommand{\bQ}{\mathbb Q}
\newcommand{\bZ}{\mathbb Z}
\newcommand{\Del}{\mathbf \Delta}
\newcommand{\var}{\varepsilon}
\newcommand{\bz}{\mathbf 0}
\newcommand{\bo}{\mathbf 1}
\newcommand{\bt}{\mathbf 2}
\newcommand{\bn}{\mathbf n}
%
%
\newcommand{\pis}{\pi_{\ast}}
\newcommand{\hpi}{\hat{\pi}_{1}}
\newcommand{\pinat}[2]{\pi^{\#}_{#1}({#2})}
\newcommand{\piA}{\pi_{\A}}
\newcommand{\coH}[3]{H^{#1}({#2};{#3})}
\newcommand{\HAQ}[3]{H^{#1}\sb{\operatorname{AQ}}({#2};{#3})}
\newcommand{\HBW}[3]{H^{#1}\sb{\operatorname{BW}}({#2};{#3})}
\newcommand{\HSO}[3]{H^{#1}\sb{\operatorname{SO}}({#2};{#3})}
%
%
\newcommand{\hG}{\hat{\Gamma}}
\newcommand{\bS}[1]{{\mathbf S}^{#1}}
\newcommand{\hX}{\hat{X}}
\newcommand{\tX}{\tilde{X}}
\newcommand{\tY}{\tilde{Y}}
\newcommand{\hY}{\hat{Y}}
%
%
\newcommand{\BL}{B\Lambda}
\newcommand{\EM}[3]{E^{#1}({#2},{#3})}
\newcommand{\EL}[2]{\EM{\Lambda}{#1}{#2}}
\newcommand{\Bd}{B_{\bullet}}
\newcommand{\Cu}{C^{\bullet}}
\newcommand{\Cuc}{C_{c}^{\bullet}}
\newcommand{\Ed}{\E_{\bullet}}
\newcommand{\Eu}{E^{\bullet}}
\newcommand{\hEd}{\hat{\E}_{\bullet}}
\newcommand{\Fu}{F^{\bullet}}
\newcommand{\Fua}{F^{\ast}}
\newcommand{\fd}{f_{\bullet}}
\newcommand{\Gd}{G_{\bullet}}
\newcommand{\gd}{g_{\bullet}}
\newcommand{\Md}{\M_{\bullet}}
\newcommand{\tNd}{\tN_{\bullet}}
\newcommand{\Qd}{Q_{\bullet}}
\newcommand{\Vd}{V_{\bullet}}
\newcommand{\tVd}{\tV_{\bullet}}
\newcommand{\Wd}{W_{\bullet}}
\newcommand{\Xd}{X_{\bullet}}
\newcommand{\tXd}{\tX_{\bullet}}
\newcommand{\Yd}{Y_{\bullet}}
\newcommand{\tYd}{\tY_{\bullet}}
\newcommand{\Zd}{Z_{\bullet}}
%
%
\newcommand{\dm}{\partial_{\max}}
\newcommand{\td}{\tilde{d}}
\newcommand{\tk}{\tilde{k}}
\newcommand{\tp}{\tilde{p}}
\newcommand{\tphi}{\tilde{\phi}}
\newcommand{\tpsi}{\tilde{\psi}}
%
%
\newcommand{\An}[1]{\A\langle{#1}\rangle}
\newcommand{\fXn}[1]{\fX\langle{#1}\rangle}
\newcommand{\Qn}[2]{Q\langle{#1}\rangle_{#2}}
\newcommand{\Qnd}[1]{\Qn{#1}{\bullet}}
\newcommand{\Xn}[1]{X\langle{#1}\rangle_{\bullet}}
\newcommand{\Xpn}[1]{\hat{X}\langle{#1}\rangle_{\bullet}}
\newcommand{\Zn}[1]{Z\langle{#1}\rangle}
\newcommand{\Po}[1]{P_{#1}}
\newcommand{\PA}[1]{P^{\A}_{#1}}
%
%
\newcommand{\fff}{\mathfrak{f}}
\newcommand{\fg}{\mathfrak{g}}
\newcommand{\fL}{\mathfrak{L}}
\newcommand{\fM}{\mathfrak{M}}
\newcommand{\fMA}{\fM_{\A}}
\newcommand{\fMAs}{\fM_{\A}^{\Stov}}
\newcommand{\fMAm}{\fM_{\A}^{\min}}
\newcommand{\fPi}{\mathbf{\Pi}}
\newcommand{\fPA}[2]{\fPi_{\A}({#1})_{#2}}
\newcommand{\fPAd}[1]{\fPA{#1}{\bullet}}
\newcommand{\fV}{\mathfrak{V}}
\newcommand{\fVd}{\fV_{\bullet}}
\newcommand{\fVn}[2]{\fV\langle{#1}\rangle_{#2}}
\newcommand{\fVnd}[1]{\fVn{#1}{\bullet}}
\newcommand{\fX}{\mathfrak{X}}
\newcommand{\fXs}{\fX^{\Stov}}
\newcommand{\fXm}{\fX^{\min}}
\newcommand{\fY}{\mathfrak{Y}}
\newcommand{\fZ}{\mathfrak{Z}}
\setcounter{section}{-1}
%
%
\title{Comparing cohomology obstructions}
\author{Hans-Joachim Baues}
\address{Max-Planck-Institut f\"{u}r Mathematik\\ 
Vivatsgasse 7\\ 53111 Bonn, Germany}
\email{baues@mpim-bonn.mpg.de}
\author{David Blanc}
\address{Department of Mathematics\\ University of Haifa\\ 31905 Haifa, Israel}
\email{blanc@math.haifa.ac.il}
\date{February 10, 2010; revised July 21, 2010}
\begin{abstract}
We show that three different kinds of cohomology \wh Baues-Wirsching
  cohomology, the \ww{\SaO}-cohomology of Dwyer-Kan, and the 
Andr\'{e}-Quillen cohomology of a \Pal\ \wh are isomorphic, under
certain assumptions. This is then used to identify the 
cohomological obstructions in three general approaches to realizability 
problems: the track category version of Baues-Wirsching, the diagram 
rectifications of Dwyer-Kan-Smith, and the \Pal\ realization of 
Dwyer-Kan-Stover. Our main tool in this identification is the notion
of a \emph{\ma}: a simplicially enriched version of an algebra over a theory.
\end{abstract}

\maketitle
%
%
\section{Introduction}

A number of questions arising in topology can be framed in terms of realizing 
an algebraic or homotopic structure in a topological setting: for
example, realizing an unstable algebra over the Steenrod algebra as
the cohomology of a space, realizing a \Pal, or lifting a group action
up to homotopy to a strict action. In these examples, the answer
appears in the form of an obstruction theory, in which elements in
appropriate cohomology groups serve both as the obstructions to
realization, and as difference obstructions which classify the various
possible realizations. 

Three general approaches to dealing with such questions have been described 
in \cite{BauCF}, \cite{DKStE,DKStB,BDGoeR}, and \cite{DKSmiR}, respectively.
Our goal in this paper is to prove that these three approaches
essentially coincide, in the cases where they all apply. In order to
do so, we introduce the notion of a \emph{\ma} \wh a
simplicially  enriched version of an algebra over a theory, in the
sense of Lawvere and Ehresmann (see Section \ref{cma}) \wh and describe a
fourth approach to the realization problem using this concept.

An important example of these methods is contained in the work of Goerss, 
Hopkins, and Miller on realizing ring spectra as structured spectra 
(cf.\ \cite{GHopkM}).

To show that the four approaches coincide, we first exhibit natural
isomorphisms between the various kinds of cohomology, after
identifying both the objects to which they apply, and the coefficient systems:

\begin{enumerate}
\renewcommand{\labelenumi}{(\alph{enumi})\ }
\item The Baues-Wirsching cohomology \w{\HBW{\ast}{\K}{D}} of a small category 
$\K$ with coefficients in a natural system $D$ (see \S \ref{sbwc}), 
\item The \ww{\SaO}-cohomology \w{\HSO{\ast}{\Z}{M}} of a simplicially
enriched category $\Z$, with coefficients in a module $M$ over the
track category  \w{\hpi{\Z}} (see \S \ref{nsacoh}).
\item The Andr\'{e}-Quillen cohomology \w{\HAQ{\ast}{\Lambda}{M}} of a 
\Pal\ $\Lambda$, with coefficients in a $\Lambda$-module $M$ (see \S
\ref{naqcoh}). 
\end{enumerate}

The identification of (a) and (b), under suitable circumstances, is given in 
Theorem \ref{tbwso}; that of (b) and (c) is given in Theorem \ref{tpaso}. 
After identifying the cohomology groups, we also identify the obstructions, 
for which we need:

\begin{mysubsection}[\label{dmodels}]{The basic setting}
Let $\C$ be a pointed model category. A collection of \emph{spherical
objects} for $\C$ is a set $\A$ of cofibrant homotopy cogroup
objects in $\C$, closed under the suspension. The motivating example
is the collection of spheres \w{\A=\{\bS{n}\}_{n=1}^{\infty}} in the
category of topological spaces, but there are many others.

Let \w{\PiA} denote the full subcategory of the homotopy category \w{\ho\C} 
whose objects are finite coproducts of objects from $\A$. A \emph{\PAa} is a 
contravariant functor \w{\Lambda:\PiA\to\Seta} which takes coproducts
to products. The category of all \PAa s is denoted by \w[.]{\PAAlg}

Such a \PAa\ $\Lambda$ is determined by its value \w{\Lambda\lin{A}\in\Seta} 
on each \w[,]{A\in\A} together with a map 
\w{\xi^{\ast}:\prod_{i\in I}\Lambda\lin{A_{i}}\to \Lambda\lin{A}} for every 
\w{\xi:A\to \coprod_{i\in I}~A_{i}} in \w[.]{\PiA\subseteq\ho\C}
Because each \w{A\in\A} is a homotopy cogroup 
object, each \w{\Lambda\lin{A}} has an underlying group structure (although
the operations \w{\xi^{\ast}} need not be group homomorphisms).

Thus when \w[,]{\A=\{\bS{n}\}_{n=1}^{\infty}} as above, a \PAa\ (called
simply a \Pal) is a graded group \w{(G_{i})_{i=1}^{\infty}} with
Whitehead products, composition operations, and a \ww{G_{1}}-action on
each \w[,]{G_{n}} as for the homotopy groups \w{\pis X} of a space $X$.

For simplicity we assume that for any collection \w{\{A_{i}\}_{i\in I}} 
of objects from $\A$ and any \w[,]{B\in\A} the natural map
\begin{myeq}[\label{eqsmall}]
\begin{split}
\colim_{J}[B,\,\coprod_{j\in J}\,A_{j}]_{\ho\C}~\longrightarrow~
[B,\,\coprod_{i\in I}\,A_{i}]_{\ho\C}
\end{split}
\end{myeq}
\noindent is an isomorphism, where the colimit on the left is taken
over the lattice of all finite subsets \w[.]{J\subseteq I} 
\end{mysubsection}

\begin{mysubsection}[\label{sthreea}]{The basic problem}
The canonical example of a \PAa\ is a \emph{realizable} one, denoted by 
\w[,]{\piA X} for fixed \w[.]{X\in\C} This is defined by setting 
\w{(\piA X)\lin{A}:=[A,X]_{\ho\C}} for each \w[.]{A\in\PiA} 

The problem we consider in this paper is that of \emph{realizing} an
abstract \PAa\ $\Lambda$: that is, finding an object \w{X\in\C} with  
\w[.]{\piA X\cong\Lambda} Such an $X$ may not exist, and need not be
unique.  There are three main approaches to the realization problem,
each describing the obstructions in terms of appropriate cohomology
classes:

\begin{enumerate}
\renewcommand{\labelenumi}{(\alph{enumi})~}
\item Trying to lift $\Lambda$ to a ``secondary \PAa'', which has 
additional structure encoding the second-order homotopy operations in
the model category $\C$ in terms of track categories. In this case,
the obstruction to such a lifting lies in Baues-Wirsching cohomology 
(see \S \ref{srelv}). 

One could try in principle to continue this process to ``higher order
track categories'', but the appropriate setting for this is not yet
clear (see \cite{BauHO} and \cite{BPaolT}).
\item Starting with a simplicial \PAa-resolution of $\Lambda$, we
obtain a ``simplicial object up to homotopy'' over $\C$. We try to
rectify it in $\C$ to a strict simplicial object. If we succeed, we
can show that its  ``geometric realization'' realizes the given \PAa\ 
$\Lambda$. 

In this setting $\Lambda$, together with \w[,]{\PiA} can be used to construct a 
certain category $\K$, as well as a simplicially enriched category,
such that the Dwyer-Kan-Smith obstructions to rectifying the
``simplicial object up to homotopy'' lie in the \ww{\SaO}-cohomology of $\K$
(see \S \ref{sdks}).
\item Starting again with a simplicial \PAa-resolution of $\Lambda$,
and trying to lift it to a strict simplicial object over $\C$  through a 
Postnikov tower, as in \cite{BDGoeR}. In this case the obstructions lie in
the Andr\'{e}-Quillen cohomology of $\Lambda$ (see Theorem \ref{treal}).
\end{enumerate}

The identification of the obstructions appearing in (a) and (b) is given in 
Theorem \ref{tzero}. In order to do this for (b) and (c), we set up
yet a fourth version of the obstruction theory in terms of \Ama s. The 
identification is then given via Theorem \ref{tobsma} and Remark
\ref{robsma}. 
\end{mysubsection}

\begin{remark}\label{rdpa}
We observe that one can dualize this setting by taking a set $\A$ of \emph{group} 
objects in \w{\ho\C} as our dual spherical objects, and define \w{\PuA} to be the 
full subcategory of \w{\ho\C} consisting of finite products of objects from $\A$. 
A \ww{\PuA}-algebra is then a covariant product-preserving functor
\w[.]{\PuA\to\Seta} This is one reason why we work in a general
categorical setting, which can readily be dualized. However,
the dual of \wref{eqsmall} is unlikely to hold, so more care is
needed in dealing with infinite products of objects from $\A$.

An important example is provided by letting \w{\A=\{K(\Fp,n)\}_{n=1}^{\infty}} 
consist of the mod $p$ Eilenberg-Mac~Lane spaces. In this case a
\ww{\PuA}-algebra is just an unstable algebra over the mod $p$
Steenrod algebra (cf.\ \cite[\S 1.4]{SchwU}). See \cite{BauAS} and
\cite{BlaCS} for more details. 
\end{remark}

\begin{mysubsection}{Organization}\label{sorg}
Section \ref{cmcso} describes the respective abstract model category
settings for the cohomology theories and the general realization problem.
Section \ref{ctrack} provides some background on track categories and the 
Baues-Wirsching cohomology of small categories. In Section \ref{csimpmc} 
we define \ww{\SaO}-categories, and show how Baues-Wirsching
cohomology can be identified with \ww{\SaO}-cohomology  (Theorem
\ref{tbwso}). In Section \ref{ccohpa} we similarly show how the
Andr\'{e}-Quillen cohomology of a \PAa\ can be identified with
relative \ww{\SaO}-cohomology (Theorem \ref{tpaso}).

In the second half of the paper, we describe the various obstruction theories 
and show how they correspond: The Dwyer-Kan-Smith \ww{\SaO}-obstructions 
to rectifying homotopy-commutative diagrams are defined in Section
\ref{cdoc}, and in Section \ref{cfobst} the  Baues-Wirsching class for
classifying linear track extensions is identified with the first
\ww{\SaO}-obstruction (Theorem \ref{tzero}). The Dwyer-Kan-Stover
approach to realizing \Pal s via Andr\'{e}-Quillen cohomology
obstructions is described in Section \ref{crpa}. In Section \ref{cma}
we introduce the concept of an \emph{\Ama}, and describe the main
example, the \Sma s, in Section \ref{csc}. Finally, Section \ref{crma}
reinterprets the obstruction theory of \cite{BDGoeR} in terms of \ma
s, and shows how they may be used to identify the Andr\'{e}-Quillen  
obstructions to realizing a \PAa\ as suitable \ww{\SaO}-obstructions 
(Theorem \ref{tobsma}).
\end{mysubsection}

\begin{mysubsection}[\label{snac}]{Notation and conventions}
The category of pointed connected topological spaces will be denoted
by \w[,]{\Ta} that of pointed sets by \w[,]{\Seta} that of
groups by \w[,]{\Gp} and that of groupoids by \w[.]{\Gpd} 
For any category $\C$, \w{s\C} denotes the category of simplicial objects 
over $\C$. However, \w{s\Set} is denoted by \w[,]{\Ss}
\w{s\Seta} by \w[,]{\Sa} and \w{s\Gp} by $\G$. The full subcategory of
\emph{reduced} simplicial sets in \w{\Sa} (with a single $0$-simplex)
will be denoted by \w[.]{\Sr} Objects in \w{s\C} will generally be
written \w[,]{\Xd} \w[,]{\Yd} and so on. The constant simplicial
object on an object \w{X\in\C} is written \w[.]{\co{X}\in s\C} 

If $\Del$ is the category of finite ordered sets \w{\bz,\bo,\bt,\dotsc} 
with order-preserving maps, then \w[.]{s\C\cong\C^{\Del}} We write
\w{\tau_{n}\Del} for the full subcategory of $\Del$ with objects
\w[,]{\{\bz,\bo,\dotsc,\bn\}} and the corresponding diagram category
\w{\C^{\tau_{n}\Del}} is called the category of $n$-\emph{truncated}
simplicial objects in $\C$, also denoted by \w[.]{\tau_{n}s\C} The
inclusion \w{\iota_{n}:\tau_{n}\Del\hra\Del} induces the
$n$-\emph{truncation} functor \w[.]{\tau_{n}:s\C\to\tau_{n}s\C} Its
left adjoint (when it exists) induces the $n$-\emph{skeleton} functor
\w[,]{\sk{n}:s\C\to s\C} and its right adjoint induces the
$n$-\emph{coskeleton} functor \w[.]{\csk{n}:s\C\to s\C} 

Given \w{A\in\Ss} and an object $X$ in a category $\C$ with coproducts, 
define \w{X\hotimes A\in s\C} by \w[,]{(X\hotimes A)_{n}:=\coprod_{a\in A_{n}} X}
with face and degeneracy maps induced from those of $A$. 
For \w[,]{Y\in s\C} set \w[.]{Y\otimes A:=\diag(Y\hotimes A)\in s\C} 

The category of all small categories will be denoted by \w[.]{\Cat}
For any set $\OO$, \w{\OC} denotes the subcategory of \w{\Cat} consisting of the 
categories having \w[,]{\Obj(\C)=\OO} with functors which are the
identity on objects.
\end{mysubsection}

\begin{ack}
We thank the referee for his or her comments.
The second author would like to thank the Max-Planck-Institut f\"{u}r
Mathematik for its repeated hospitality while this research was
carried out.  He would also like to thank Bernard Badzioch, Wojtek
Dorabia{\l}a, Mark Johnson, and Jim Turner for many useful
discussions.
\end{ack}

%
%
\sect{Model categories and cohomology}
\label{cmcso}

We first describe the model category framework needed to define the 
cohomology theories, and study the realization problems described above:

\begin{assume}\label{amodcat}
We assume throughout this paper that our model categories are pointed, 
cofibrantly generated, simplicial (see \cite[II, \S 1]{QuiH}), and right proper 
(that is, the pullback of a weak equivalence along a fibration is a
weak equivalence).  

For simplicity of treatment we will assume that all objects in $\C$ are fibrant
(although many of our constructions make use of the category $\Ss$ of
simplicial sets,  where this does not hold). Note that we may take
\w{\C=\G} if we want such a model category for the homotopy theory of
pointed connected topological spaces. 
\end{assume}

First, in order to provide an appropriate setting for resolutions, we
shall need to deal with simplicial objects over our model category
$\C$, for which we have the following:

\begin{defn}\label{drmc}
Let $\C$ be a model category as above, with a set $\A$ of spherical
objects (\S \ref{dmodels}). In order to define a model category
structure on \w[,]{s\C} we choose the set 
\w{\tA:=\{A\otimes\Del[n]/(A\otimes\partial\Del[n])\}_{n\in\bN,A\in\A}} 
(see \S \ref{snac}) as the collection of spherical objects for \w[.]{s\C} 
We think of \w{A\otimes\Del[n]/(A\otimes\partial\Del[n])} as
the \emph{simplicial} suspension of $A$; we reserve the notation
\w{\Sigma^{k}} for (internal) suspension in $\C$.   

Extending the simplicial structure from $\C$ to \w{s\C} in the usual
way (cf.\ \cite[II, \S 4]{QuiH}), we set 
\w{[\Xd,\Yd)_{\ho s\C}:=\pi_{0}\map_{s\C}(\Xd,\Yd)} for
\w[.]{\Xd,\Yd\in s\C} We write \w{\pinat{n}{\Xd}} for the $\A$-graded
group \w{[A\hotimes\bS{n},\Xd]_{\ho s\C}} \wb[.]{A\in\A} These are
called the \emph{natural}  homotopy groups of \w[.]{\Xd}

We now define the \emph{resolution model category structure} on \w{s\C}
determined by $\A$, by letting a simplicial map \w{f:\Xd\to\Yd} be:
\begin{enumerate}
\renewcommand{\labelenumi}{(\roman{enumi})}
\item a \emph{weak equivalence} if \w{\pinat{\ast}{f}} is a weak
equivalence of $\A$-graded simplicial groups.
\item a \emph{cofibration} if it is (a retract of) a map with the
following property: for each \w[,]{n\geq 0} there is a cofibrant 
object \w{W_{n}} in $\C$  which is weakly equivalent to a coproduct
of objects from $\A$,  and a map  \w{\varphi_{n}:W_{n}\to Y_{n}} in
$\C$ inducing a trivial cofibration
\w[.]{(X_{n}\amalg_{L_{n}\Xd}L_{n}\Yd)\amalg W_{n}\to Y_{n}} Here
\w{L_{n}\Yd} is the $n$-th latching object for \w{\Yd} 
(cf.\ \cite[\S 2.1]{BJTurR}).
\item a \emph{fibration} if it is a Reedy fibration (cf.\
\cite[15.3]{PHirM}) and \w{\piA f} (\S \ref{sthreea}) is a fibration
of $\A$-graded simplicial groups. 
\end{enumerate}
\noindent See \cite{BousCR} and \cite{DKStE}.
\end{defn}

Applying \w{\piA} in each simplicial dimension to any \w{\Xd\in s\C}
yields a simplicial \PAa\ \w[.]{\piA\Xd} By taking the usual homotopy
groups of the underlying $\A$-graded simplicial group in each degree,
we obtain the $\bN$-graded \PAa\ \w[.]{\pis\piA\Xd} This is related to
natural homotopy groups by a \emph{spiral long exact sequence} (cf.\
\cite[8.1]{DKStB}): 
\begin{myeq}[\label{eqspiral}]
\begin{split}
\ldots~\to~\Omega\pinat{n-1}{\Xd}&~\xra{s_{n}}~
\pinat{n}{\Xd}~\xra{h_{n}}~\pi_{n}\piA\Xd~\xra{\partial_{n}} 
\Omega\pinat{n-2}{\Xd}\\ ~& \xra{s_{n-1}}~\pinat{n-1}{\Xd}
\to~\ldots~\to~\pinat{0}{\Xd}~\xra{\cong}~\pi_{0}\piA\Xd 
\end{split}
\end{myeq}

It follows that a map \w{f:\Xd\to\Yd} in \w{s\C} is a weak equivalence if 
and only if the map of simplicial \PAa s \w{f_{\ast}:\piA\Xd\to\piA\Yd} is a weak
equivalence in the resolution model category \w[.]{s\PAAlg}

\begin{mysubsection}{Examples of resolution model categories}\label{egrmc}

\begin{enumerate}
\renewcommand{\labelenumi}{(\alph{enumi})\ }
\item Let \w{\C=\Gp} with the \emph{trivial} model category structure, and 
\w[.]{\A:=\{\bZ\}} The resulting resolution model category structure on 
\w{\G:= s\C} is the usual one.
\item More generally, let \w{\C=\TAlg} be a category of universal algebras
(with an underlying group structure), represented by a theory $\Theta$ 
(cf.\ \cite[\S 1]{ARosiL}), such as \w[.]{\PAAlg} In this case we let $\A$
be the collection of free monogenic algebras.
\item We can iterate the process by taking $\G$ for $\C$, and letting $\A$
consist of the $\G$-spheres. We thus obtain a resolution model category structure 
on \w{s\G} (or on \w[),]{s\Ta} which is the original example of \cite{DKStE}.
\item If $\C$ is a resolution model category and $I$ is some small
category, the category \w{\C^{I}} of $I$-diagrams in $\C$ also has
a resolution model category structure, in which the spherical objects
are certain free $I$-diagrams (cf.\ \cite[\S 1]{BJTurR}).
\end{enumerate}
\end{mysubsection}

In order to define cohomology groups in our model category, it is 
convenient to consider the following setting:

\begin{defn}\label{dssmc}
A model category $\C$ is called \emph{semi-spherical} (see \cite[\S 2.23]{BJTurH})
if it is equipped with:
\begin{enumerate}
\renewcommand{\labelenumi}{(\alph{enumi})\ }
\item A \emph{coefficient category} \w[,]{\Coef(\C)} together a functor
\w[.]{\hpi:\C\to\Coef(C)}
\item For each \w[,]{n\geq 2} a functor \w{\pi_{n}:\C\to\RM{\hpi(-)}} taking 
\w{Z\in\C} into the category of \emph{modules} over \w{\hpi Z} (that is,
abelian group objects in \w[).]{\Coef(C)/\hpi Z}
\item Each \w{Z\in\C} has a functorial \emph{Postnikov tower} of
  fibrations under $Z$: 
$$
Z \dotsc \to \Po{n}Z\xra{p\q{n}}\Po{n-1}Z\xra{p\q{n-1}}\dots\to\Po{0}Z~,
$$
\noindent with \w{Z\to\lim_{n}\Po{n}Z} a weak equivalence, and the usual
properties for the structure maps \w[.]{r\q{n}:Z\to\Po{n}Z}
\item For every \w[,]{\Lambda\in\Coef(\C)} there is a functorial
  \emph{classifying object} \w[,]{\BL\in\C} unique up to homotopy, with 
  \w{\BL\simeq\Po{1}\BL} and \w[.]{\hpi\BL\cong \Lambda}
\item Given \w{\Lambda\in\Coef(\C)} and a $\Lambda$-module $M$, for
  each \w{n\geq 1} there is a functorial \emph{Eilenberg-Mac~Lane
  object} \w{E=\EL{G}{n}} in $\C$, unique up to homotopy, equipped
  with a section $s$ for \w[,]{r\q{1}:E\to\Po{1}E\simeq \BL} such that
  \w{\pi_{n}E\cong M} as $\Lambda$-modules and \w{\pi_{k}E=0} for
  \w[.]{k\neq 0,1,n} 
\item For every \w[,]{n\geq 1} there is a functor that assigns to each 
\w{Z\in\C} a homotopy pull-back square:
%
\mydiagram[\label{eqkinv}]{
\ar @{} [dr] |<<<{\framebox{\scriptsize{PB}}}
\Po{n+1}Z \ar[r]^{p\q{n+1}} \ar[d] &
\Po{n}Z \ar[d]^{k_{n}}\\ B(\hpi Z) \ar[r] & \EM{\hpi Z}{\pi_{n+1}Z}{n+2}~.
}
\noindent The map \w{k_{n}} is called the $n$-th $k$-\emph{invariant} for $Z$.
\end{enumerate}
\end{defn}

\begin{examples}\label{egssmc}
The motivating example is the category \w{\Ta} of pointed topological spaces.

In addition, all the resolution model categories of \S \ref{egrmc} are
semi-spherical (see \cite[\S 3]{BJTurR}). We note that for the ``algebraic''
categories \w{\C=s\TAlg} of simplicial universal algebras (\S \ref{egrmc}(b)),
\w{\hpi{\Xd}} is just \w[,]{\pi_{0}\Xd}and \w{\Coef(\C)} is \w{\TAlg} itself.
Thus a module over a \Tal\ $\Lambda$ is just an abelian group object in 
\w{\TAlg/\Lambda} (cf.\ \cite{BecT}).
\end{examples}

\begin{defn}\label{dcoh}
Let $\C$ be a semi-spherical simplicial model category, and assume given 
\w[,]{\Lambda\in\Coef(\C)} a $\Lambda$-module $M$, and  an object 
\w{Z\in\C}  equipped with a \emph{twisting map}  \w[.]{p:\hpi Z\to \Lambda} 
Following \cite[II, \S 5]{QuiH}, we define the $n$-th \emph{cohomology group} 
of $Z$ with coefficients in $M$ to be
$$
\coH{n}{Z/\Lambda}{M}~:=~[Z,\EL{M}{n}]_{\C/\BL}~=~
\pi_{0}\map_{\C/\BL}(Z,\EL{M}{n})~,
$$
\noindent where \w{\map_{\C/A}(Z,Y)} is the sub-simplicial set of 
the mapping space \w{\map_{\C}(Z,Y)} in $\C$ consisting of maps over a
fixed base $A$. 

Typically, we have \w[,]{\Lambda=\hpi Z} with $p$ a weak equivalence; if 
in addition \w[,]{Z\simeq \BL} we denote \w{\coH{n}{Z/\Lambda}{M}} simply by 
\w[.]{\coH{n}{\Lambda}{M}} 
\end{defn}

\begin{remark}\label{rrelcoh}
There is also a relative version, for a cofibration \w{i:X\hra Y} in \w[:]{\C/\BL}
If $Z$ is the cofiber of $i$ in \w{\C/\BL} \wh that is, the homotopy 
pushout of:
\mydiagram[\label{eqpo}]{
\ar @{}[drr] |>>>>>{\framebox{\scriptsize{PO}}}
X \ar[d]_{p} \ar[rr]^{i} & & Y \ar[d]^{q}\\
\BL \ar[rr] & & Z~,
}
\noindent then
$$
\coH{n}{(Y,X)/\Lambda}{M}~:=~[(Z,\BL),\,(\EL{M}{n},\BL)]_{\C/\BL}~.
$$
\noindent (cf.\ \cite[\S 2.1]{DKSmiO}). Again if \w{\Lambda=\hpi Y} we
write simply \w[.]{\coH{n}{Y,X}{M}}
\end{remark}

\begin{mysubsection}{The module \ww{\Omega\Lambda}}\label{emodule}
We close this section with the following important example of a module in the
category of \PAa s: 

Given a \PAa\ $\Lambda$, we define the \PAa\ \w{\Omega_{+}\Lambda} as
an $\A$-graded group by
\w[.]{(\Omega_{+}\Lambda)\lin{A}:=\Lambda\lin{\Sigma A\vee A}} 
We identify the \PAa\ structure on \w{\Omega_{+}\Lambda}  as follows:

Given \w{f:B\to A} in \w[,]{\PiA} define \w{\nabla f:B\to A\vee A} to be
\w[,]{-i_{2}\circ f +(i_{1}+i_{2})\circ f} using the
co-group structure on $B$ (where \w{i_{1},i_{2}:A\to A\vee A} are the
two inclusions). If \w{j:A\vee A\hra A} is \w[,]{\ast\vee\Id} then 
\w[,]{j\circ\nabla f\sim\ast} with a nullhomotopy \w[.]{H:CB\to A} 
Now let \w{I_{\ast}X} denote the reduced cylinder in $\C$
and let $G$ be the composite of
$$
I_{\ast}B~\xra{I_{\ast}\nabla f}~I_{\ast}(A\vee A)=I_{\ast}A\vee I_{\ast}A~
\xra{q\vee p_{0}}~\Sigma A\vee A~,
$$
\noindent where \w{q:I_{\ast}A\to\Sigma A} is the quotient map and
\w{p_{0}:I_{\ast}A\to A} is the projection. If we identify \w{\Sigma
  B} with the pushout \w{CB\cup_{B}I_{\ast}B\cup_{B}CB} (under the two
inclusions of $B$ into \w[),]{I_{\ast}B} we define 
\w{E(\nabla f):\Sigma B\to \Sigma A\vee A} to be the map given on the
pushout by \w[,]{(H,G,H)} and call it the \emph{partial suspension} of
\w[.]{\nabla f}  Because \w{[\Sigma B,\Sigma A\vee A]} is an abelian
group, this  is independent of $H$. See \cite[\S 3]{BauO} for more
details, including explicit rules for applying the partial suspension
to maps among wedges of spheres. 

Since $\Lambda$ is contravariant, the map 
\w{(E\nabla f, i_{2}\circ f):\Sigma B\vee B\to\Sigma A\vee A}  induces
the required map
\w[.]{f^{\ast}:\Omega_{+}\Lambda\lin{A}\to\Omega_{+}\Lambda\lin{B}}  
We thus have a split exact sequence of \PAa s:
$$
\xymatrix@R=25pt{
\ast \rto & \Omega\Lambda \rto & \Omega_{+}\Lambda \rto_{p} &
~\ar@/_1pc/[l] \Lambda \rto & \ast~,}
$$
\noindent where \w[.]{\Omega\Lambda:=\Ker(p)} This gives
\w{\Omega_{+}\Lambda} the structure of a module over $\Lambda$ \wh or
equivalently, a natural system on \w{\PiA} (see \S \ref{dnatsys}
below). Note that \w[.]{(\Omega\Lambda)\lin{A}\cong\Lambda\lin{\Sigma A}} 
for all \w[,]{A\in\PiA} but the operation
\w{f^{\ast}:\Omega\Lambda\lin{A}\to\Omega\Lambda\lin{B}} described as
above for \w[,]{f:B\to A} is \emph{not} in general 
\w[.]{(\Sigma f)^{\ast}:\Lambda\lin{\Sigma A}\to\Lambda\lin{\Sigma B}} 
\end{mysubsection}

\begin{remark}\label{rtwoloop}
A canonical identification of \w{(A\otimes\bS{1})/(\ast\otimes\bS{1})} with 
\w{\Sigma A\vee A} in any pointed model category is given in
\cite{BJiblSL}, such that:
$$
\xymatrix@R=25pt{
(B\otimes\bS{1})/(\ast\otimes\bS{1})
\ar[rrr]^{(f\otimes\bS{1})/(\ast\otimes\bS{1})} 
\ar[d]_{\simeq} &&& (A\otimes\bS{1})/(\ast\otimes\bS{1}) \ar[d]_{\simeq}\\
\Sigma B\vee B \ar[rrr]_{(E\nabla f, i_{2}\circ f)} &&& \Sigma A\vee A
}
$$
\noindent commutes up to homotopy for any \w[.]{f:B\to A} Thus our
definition of \w{\Omega_{+}\Lambda} agrees with that of \cite[\S 9.4]{DKStB}.
\end{remark}

%
%
\sect{Track categories and natural systems}
\label{ctrack}

The first approach to realization problems in \S \ref{sthreea}(a), 
developed in \cite{BWirC} (cf.\ \cite[\S 2-3]{BauAS}), concentrates on 
secondary homotopy structure, in the following sense:

\begin{defn}\label{dtrackcat}
A \emph{track category} is a category $\E$ enriched in groupoids.
It thus consists of two categories \w{\E_{0}} and \w[,]{\E_{1}} with the same 
objects, and two functors \w{s,t:\E_{1}\to\E_{0}} which are the identity
on objects. Here \w{\E_{0}} is the ordinary category underlying $\E$, 
while \w{\E_{1}(X,Y)} is a groupoid, with maps
\w{f:X\to Y} in \w{\E_{0}} as objects, and a set \w{\E_{1}(f,g)} of 
morphisms (called \emph{$2$-cells}) from \w{s(H)=f:X\to Y} to 
\w[,]{t(H)=g:X\to Y} written \w[.]{H:f\Ra g} The groupoid operation 
is denoted by \w[,]{H\Box H'} when defined.

There is natural equivalence relation on maps in \w{\E_{0}} induced by 
the $2$-cells, and the quotient category \w{\ho\E} is called the 
\emph{homotopy category} of $\E$. See \cite[VI, \S 3]{BauCHF} for
further details. 
\end{defn}

\begin{example}\label{egtrackcat}
The motivating example of a track category is obtained from a model 
category $\C$ by letting \w{\E_{0}=\C_{\cf}} (the full subcategory of 
fibrant and cofibrant objects), with \w{\E_{1}(X,Y)} the groupoid of tracks
(homotopy classes of homotopies) between maps from $X$ to $Y$ in $\C$. This is 
called the \emph{homotopy track category} of $\C$.
\end{example}

\begin{remark}\label{rtrackcat}
There is a model category structure on the category \w{\Track} 
of (small) track categories, with (strict) track functors, in which 
the weak equivalences are bi-essential surjections \w{F:\E\to\E'} 
which induced equivalences of categories 
\w{F:\E_{1}(X,Y)\xra{\cong}\E'_{1}(FX,FY)} (see \cite{LackQ}).
\end{remark}

\begin{defn}\label{dnatsys}
For any category $\K$, the \emph{category of factorizations} of $\K$ is
the category \w{\Fac\K} having as objects \w{\Arr\K} (the morphism set
of $\K$) and as morphisms from $f$ to $g$ commuting squares of the form:
$$
\xymatrix@R=25pt{
1 \ar[r]^{f} & 2 \ar[d]^{\beta}\\
0 \ar[u]_{\alpha} \ar[r]_{g} & 3
}
$$
\noindent with the obvious composition.  A \emph{natural system on $\K$} with 
values in a category $\M$ is a functor \w[.]{D:\Fac\K\to\M} The category 
of such natural systems (with natural transformations as morphisms) will 
be denoted by \w[.]{\NS_{\K}(\M)=\M^{\Fac\K}} When \w[,]{\M=\Abgp} $D$ is 
called simply a \emph{natural system on $\K$} (see \cite[\S 1]{BWirC}).
\end{defn}

\begin{example}\label{egnatsys}
For \w{\A\subseteq\C} as in \S \ref{dmodels}, we have a canonical natural
system \w{\Omega\PiA} on \w[,]{\PiA} defined for \w{f:B\to A} in \w{\PiA} by
\w[.]{\Omega\PiA(f):=\Hom_{\PiA}(\Sigma B, A)=[\Sigma B, A]_{\ho\C}}
For \w[,]{g:A\to A'} the induced map \w{g_{\ast}:\Omega\PiA(f)\to\Omega\PiA(gf)}
is given by post-composition, while for \w[,]{h:B'\to B} 
\w{h^{\ast}:\Omega\PiA(f)\to\Omega\PiA(fh)} is given by
\w{(E\nabla h)^{\ast}(\alpha,f)} (cf.\ \S \ref{emodule}).
\end{example}

\begin{mysubsection}{Baues-Wirsching cohomology of a small category}\label{sbwc}
If \w{N\K} is the nerve of $\K$, define \w{\dm:N_{n}\K\to N_{1}\K} by
\w{\dm(\sigma):=d_{1}d_{2}\dotsb d_{n-1}\sigma\in N_{1}\K} (the composite 
of the corresponding composable sequence in $\K$), and set
$$
N_{n}[f]~:=~\{\sigma\in N_{n}\K~:\ \dm(\sigma)=f\}\hsp\text{for any arrow}\hsm 
f~\text{in}~\K~.
$$
\noindent This defines a collection of sets indexed by \w[.]{\Arr\K} 
Note that there is a forgetful functor from natural systems on $\K$ to 
\ww{\Arr\K}-graded sets, whose left adjoint is the \emph{free natural
system} functor (cf.\ \cite[\S 5.14]{BPaolT}. Thus for each \w{n\geq 0} 
we have a free natural system in sets on $\K$ denoted by \w[.]{\tN_{n}\K} 

The face and degeneracy maps of \w{N\K} induce maps of natural systems 
as follows:
\begin{enumerate}
\renewcommand{\labelenumi}{(\alph{enumi})\ }
\item If \w{\phi=d_{i}:N_{n}\K\to N_{n-1}\K} \wb{0<i<n} or
\w{\phi=s_{j}:N_{n}\K\to N_{n+1}\K} \wb[,]{0\leq j\leq n}
we define \w{\tphi:\tN_{n}\K\to \tN_{n\pm 1}\K} to be \w[.]{\F\phi}
\item Given \w[,]{\sigma\in N_{n}\K} define the map of natural systems 
\w{\td_{0}:\tN_{n}\K\to\tN_{n-1}\K} by setting 
\w[.]{\td_{0}(\sigma):=(d_{2}\dotsb d_{n}\sigma)^{\ast}(d_{0}\sigma)}
This extends to all of \w{\tN_{n}\K} by the adjointness of $\U$ and $\F$ above.
\item We similarly define the $n$-th face map \w{\td_{n}:\tN_{n}\K\to\tN_{n-1}\K} by 
\w{\td_{n}(\sigma):=(d_{0}\dotsb d_{n-2}\sigma)_{\ast}(d_{n}\sigma)}.
\end{enumerate}
This makes \w{\tNd\K:=(\tN_{n}\K)_{n=0}^{\infty}} into a simplicial
object in the category \w[.]{\NS_{\K}(\Set)}

Finally, a natural system (in \w[)]{\Abgp} on $\K$ can be thought of as
an abelian group object $\D$ in \w[,]{\NS_{\K}(\Set)} so we can define 
a cosimplicial abelian group \w{\Cu(\K;D)} by setting 
\w[.]{C^{n}(\K;D):=\Hom_{\NS_{\K}(\Set)}(\tN_{n}\K,D)} Its $n$-th
cohomotopy group is defined to be the $n$-th \emph{Baues-Wirsching
cohomology group} of $K$ with coefficients in $D$, written
\w[.]{\HBW{n}{\K}{D}~:=~\pi^{n}(\Cu(\K;D))} 

The cochain complex \w{\Fua(\K,D)} used in \cite{BWirC} to define 
\w{\HBW{n}{\K}{D}}is that associated to \w[,]{\Cu(\K;D)} so 
\w[.]{\HBW{n}{\K}{D}\cong H^{n}\Fua(\K,D)}
\end{mysubsection}

\begin{defn}\label{dlintrack}
A \emph{linear track extension} of a category $\K$ by a natural system $D$
is a track category $\E$ with \w[,]{\ho\E=\K} for which \w{\Aut_{\E}(f)}
is naturally isomorphic to \w{D([f])} for all maps $f$ in \w[.]{\E_{0}}
Such an extension is denoted by \w[.]{D\to\E\to\K}
\end{defn}

\begin{prop}[\protect{\cite[VI, Theorem 3.15]{BauCHF}}]\label{tlintrext}
The set of all linear track extensions of a category $\K$ by a given 
natural system $D$, up to ($D$-equivariant) weak equivalence, is in 
one-to-one correspondence with \w[.]{\HBW{3}{\K}{D}}
\end{prop}

This can be interpreted as describing the homotopy equivalence classes 
in \w{\Track/\ho\E} (as in \S \ref{rtrackcat}).

\begin{remark}\label{rcocycle}
If $\A$ is a set of spherical objects in a model category $\C$, 
let \w{\tCA} be a sub-track category of the homotopy track category of 
$\C$ with \w[.]{\ho\tCA\cong\PiA} This is a linear track extension 
\w[,]{\Omega\PiA\to\tCA\to\PiA} and one can describe an explicit cocycle 
representing the corresponding cohomology class \w{\chi_{\tCA}} in 
\w{\HBW{3}{\PiA}{\Omega\PiA}} as follows:

Choose an arbitrary fixed representative \w{s\phi:0\to 1} in \w{\C}
for each $1$-simplex \w{\phi:0\to1} in \w[,]{\N_{1}(\PiA)} and a fixed
track \w{H_{(\phi,\psi)}:s\phi\circ s\psi\simeq s(\phi\psi)} for each
$2$-simplex \w{0\xra{\psi}1\xra{\phi}2} in \w[.]{\N_{2}(\PiA)} Now we
associate to each $3$-simplex
\w{0\xra{\phi_{1}}1\xra{\phi_{2}}2\xra{\phi_{3}}3} in \w{\N_{3}(\PiA)} 
the element 
\begin{myeq}\label{eqbauescocyc}
H_{(\phi_{3},\phi_{2}\circ\phi_{1})}\Box
(\phi_{3})_{\ast}H_{(\phi_{2},\phi_{1})}\Box
(\phi_{1})^{\ast}H_{(\phi_{3},\phi_{2})}^{-1}\Box
H_{(\phi_{3}\circ\phi_{2},\phi_{1})}^{-1}
\end{myeq}
\noindent in \w[.]{\Aut(s(\phi_{3}\circ\phi_{2}\circ\phi_{1}))~\cong~
(\Omega\PiA)(\phi_{3}\circ\phi_{2}\circ\phi_{1})} 
\end{remark}

%
%
\sect{$\SaO$-categories and \ww{\SaO}-cohomology}
\label{csimpmc}

For the second approach to the realization problem of \S \ref{sthreea},  
due to Dwyer and Kan, we use the framework of simplicially enriched categories:

\begin{defn}\label{dsv}  
For a fixed set $\OO$, a category $\Z$ enriched in simplicial sets with object 
set $\OO$  will be called an \ww{\SO}-\emph{category}, and the category of 
all such will be denoted by \w[.]{\SOC} Equivalently, such a category $\Z$ can 
be thought of as a simplicial object in \w{\OC}  (\S \ref{snac}): this means 
$\C$ has a fixed object set $\OO$  in each dimension, and all face and degeneracy 
functors the identity on objects. 

More generally, if \w{(\V,\otimes)} is any monoidal category, a
\ww{\VO}-\emph{category} is a small category \w{\C\in\OC} enriched
over $\V$. The category of all such categories will be denoted by
\w[.]{\VOC} Examples for \w{(\V,\otimes)} include \w[,]{\TT}
\w[,]{\Gp} \w[,]{\Gpd} and $\Ss$, with \w{\otimes=\times} (Cartesian
product), or the category \w{\Set^{\Box}} of \emph{cubical sets} with
its monoidal enrichment $\otimes$ (see \cite[\S 1.5]{BJTurH}).

The main example we shall be working with is \w[,]{\V=\Sa} with 
\w{\otimes=\wedge} (smash product). Again we can identify an
\ww{\SaO}-category with a simplicial pointed $\OO$-category.
\end{defn}

\begin{mysubsection}{$\SaO$-categories}\label{ssoc}
In \cite[\S 1]{DKanS}, Dwyer and Kan define a simplicial model
category structure on \w[,]{\SOC} also valid for \w{\SaOC} (cf.\
\cite[Prop.~1.1.8]{HovM}), in which a map \w{f:\X\to\Y} is a fibration
(respectively, a weak equivalence) if for each \w[,]{a,b\in\OO} the
induced map \w{f_{(a,b)}:\X(a,b)\to\Y(a,b)} is such. 

The cofibrations in \w{\SOC} or \w{\SaOC} are not easy to describe. 
However, if \w{\K\in\OC} is any category with object set $\OO$, then 
\w{\co{\K}\in s\OC\cong\SOC} has a cofibrant replacement defined as follows:

There is a forgetful functor \w{U:\Cat\to\DG} to the category of directed 
graphs, whose left adjoint \w{F:\DG\to\Cat} is the free category
functor (cf.\ \cite[\S 2.4]{DKanS} and \cite[\S 2]{CPorV}).  Both $U$
and $F$ are the identity on objects. A canonical cofibrant replacement
for the constant simplicial category \w{\co{\K}\in s\OC} is 
provided by the simplicial category \w[,]{\Fs\K} obtained by iterating
the comonad \w{FU:\OC\to\OC} (so \w[).]{\F_{n}\K:=(FU)^{n+1}\K} 
The augmentation \w{\Fs\K\to\K} induces a weak equivalence 
\w{\Fs\K\simeq\co{\K}} in \w[.]{s\OC\approx\SOC} If $\K$ is pointed, 
\w{\Fs\K} is an \ww{\SaO}-category.

Both \w{\SOC} and \w{\SaOC} are semi-spherical (\S \ref{dssmc}), with
coefficient category \w{\GOC} (=track categories with object set $\OO$). 

The \emph{fundamental track category} of a (fibrant)  \ww{\SO}- or 
\ww{\SaO}-category $\Z$ is obtained by applying the fundamental groupoid 
functor \w{\hpi:\Ss\to\Gpd} to each mapping space \w[,]{\Z(a,b)} noting that 
\w{\hpi} commutes with cartesian products, and thus extends to \w{\SOC} (and to
\w[,]{\SaOC} too, since in the pointed case the composition factors
through $\wedge$). For each \w{n\geq 2} we obtain a \ww{\hpi Z}-module
by applying \w{\pi_{n}(-)} to each mapping space of $Z$ (again, \w{\pi_{n}}
preserves products). 

The usual Postnikov tower functor, classifying space, and Eilenberg-Mac~Lane 
functors for $\Ss$ or \w{\Sa} similarly preserve products, and thus extend 
to \w{\SOC} and \w[.]{\SaOC}  For the functorial $k$-invariants, use the
construction of \cite[\S 6]{BDGoeR}.
\end{mysubsection}

\begin{notn}\label{nsacoh}
We write \w{\HSO{\ast}{\Z/\Lambda}{M}} (or just
\w[)]{\HSO{\ast}{\Z}{M}} for the cohomology groups of an
\ww{\SaO}-category $\Z$, as defined in \S \ref{dcoh}. Similarly, we
write \w{\HSO{\ast}{(\Z,\Y)/\Lambda}{M}} (or just
\w[)]{\HSO{\ast}{\Z,\Y}{M}} for the relative cohomology of \S
\ref{rrelcoh}. We call this \ww{\SaO}-\emph{cohomology} (compare
\cite{DKSmiO}). 
\end{notn}

\begin{defn}\label{dwk}
A cubical version of the free simplicial category \w{\Fs\K} on a
category \w{\X\in\OC} is provided by the \emph{bar construction} of
Boardman and Vogt: this is a category \w{W\K} enriched in
the monoidal category \w{(\Set^{\Box},\otimes)} of cubical sets (\S \ref{snac}). 
For \w[,]{a,b\in\OO=\Obj\K} the cubical mapping complex \w{W\K(a_{n+1},a_{0})} 
has an $n$-cube \w{I^{n}\li{\fd}} for each sequence: 
%
\begin{myeq}\label{eqcompseq}
\fd~=~\left(a_{n+1}\xra{f_{n+1}}a_{n}\xra{f_{n}}a_{n-1}\dotsc
a_{1}\xra{f_{1}}a_{0}\right)~.
\end{myeq}
\noindent of \w{(n+1)} composable maps in $\K$.

The $i$-th $0$-face \w{d_{i}^{0}} of \w{I^{n}{\fd}} is identified with
  \w[,]{I^{n-1}{f_{1}\circ\dotsc\circ (f_{i}\cdot f_{i+1})\circ\dotsc f_{n+1}}}
  that is, we carry out (in $\K$) the $i$-th composition in the
  sequence \w[.]{\fd}

The cubical composition
$$
W\K(a_{0},a_{i})\otimes W\K(a_{i},a_{n+1})~\to~W\K(a_{0},a_{n+1})=W\K(a,b)
$$
identifies the ``product'' \ww{(n-1)}-cube
\w{I^{i}{f_{0}\circ\dotsc\circ f_{i}}\otimes
       I^{n-i-1}{f_{i+1}\circ\dotsc\circ f_{n+1}}}
with the $i$-th $1$-face \w{d_{i}^{1}} of \w[.]{I^{n}{\fd}}
See \cite[III, \S 1]{BVogHI} or \cite[\S 3.1]{BJTurH} for further details. 
\end{defn}

\begin{lemma}\label{lwkfs}
For any small category $\K$, the simplicial category \w{\Fs\K} is a natural
triangulation of \w[.]{W\K} 
\end{lemma}

\begin{proof}
The $n$-cube \w{I^{n}{\fd}} is subdivided into \w{n!} $n$-simplices by
fully parenthesizing \w{(f_{1},\dotsc,f_{n+1})} in all possible ways,
with the $i$-th face map defined by omitting the $i$-th level of
parentheses (cf.\ \cite[\S 2.21]{BMarkH}). 
\end{proof}

\begin{example}\label{egwkfs}
For \w[,]{n=2} given three composable maps \w[,]{0\xra{h}1\xra{g}2\xra{f}3}  
we have: 
\mydiagram[\label{eqsquare}]{
(f)(g)(h)\bullet \ar[rrr]^{((f)(g))((h))} \ar@<3.6ex>[ddd]_{((f))((g)(h))} 
\ar[rrrddd]_<<<<<<<<<<<<<<<<<<<<<<<<<{((f)(g)(h))} &&& 
\bullet(fg)(h) \ar@<-3.6ex>[ddd]^{((fg)(h))} \\
&& *[F]{~(((f)(g))((h)))~} & \\
&  *[F]{~(((f)((g)(h)))~} && \\
(f)(gh) \bullet \ar[rrr]_{((f)(gh))} &&& \bullet(fgh)
}
\end{example}

\begin{remark}\label{rnatmod}
If $D$ is a natural system on a category $\K$, with \w[,]{\OO=\Obj(\K)}
it can be thought of as an abelian group object on \w[.]{\OC/\K} Moreover,
$\K$ itself is the (discrete) fundamental groupoid of the homotopically trivial
simplicial category \w{\Fs\K\simeq\K} in \w[.]{\SOC} (or \w[,]{\SaOC} if $\K$ 
is pointed). Thus $D$ is just a module $M$ over $\K$.
\end{remark}

%
%
\begin{thm}\label{tbwso}
If $D$ is a natural system on a small pointed category $\K$, the 
$n$-th Baues-Wirsching cohomology group \w{\HBW{n}{\K}{D}} is naturally 
isomorphic to the \ww{(n-1)}-st \ww{\SaO}-cohomology group 
\w[,]{\HSO{n-1}{\K}{D}} for each \w[.]{n\geq 1}
\end{thm}

In \cite[Theorem 5.3]{DKanH}, Dwyer and Kan prove a similar result,
using a different definition of the cohomology of a small category,
which they call \emph{Hochschild-Mitchell} cohomology.

\begin{proof}
The \ww{\SaO}-cohomology groups 
\w{\HSO{n}{\K/\G}{D}\cong[\Fs\K,\,\EM{G}{M}{n}]_{\SaOC/B\G}} of \S \ref{dcoh} 
may be computed as the cohomotopy groups of the cosimplicial abelian group 
\w{\Eu:=\Hom_{\SaOC/\K}(\Fs\K,D)} (cf.\ \cite[Proposition 3.11]{BJTurHA}).

In order to compare \w{\Eu} with \w{\Cu(\K,D)} of \S \ref{sbwc}, note that
for \w[,]{n\geq 1} there is an obvious one-to-one correspondence between the
$n$-cubes of \w{W\K} (\S \ref{dwk}) and the \ww{(n+1)}-simplices of the nerve
\w[.]{\N(\K)} Moreover, for \w{n\geq 2} this extends to the face maps,
if we omit the \ww{d_{i}^{1}}-faces with \w[]{1<i<n} \wh that is,
those which are cubical products of two lower-dimensional cubes.  
There are \w{2n-(n-2)=n+2} remaining \ww{(n-1)}-facets, of which two
are the Cartesian products 
\w{I^{0}\li{f_{0}}\times I^{n}\li{f_{1},\dotsc,f_{n+1}}} and 
\w{I^{n}\li{f_{0},\dotsc,f_{n}}\times I^{0}\li{f_{n}}} (corresponding
to \w{d_{0}\sigma\li{\fd}} and \w[,]{d_{n+1}\sigma\li{\fd}}
respectively), and the others are obtaining out the adjacent
compositions as for \w{d_{i}\sigma\li{\fd}} \wb[.]{i=1,\dotsc,n}
Note that the facets we have omitted are not relevant for the
coboundary of a cubical \ww{(n-1)}-chain.

Finally, the cubical cochain complex \w{\Cuc:=\Hom_{\COC/\K}(W\K,D)} has
the same cohomology as \w{\Eu} by the Lemma \ref{lwkfs} and the
Acyclic Model Theorem (cf.\ \cite{EMacLAM}), and clearly has the same
cohomolgy as \w{\Cu(\K,D)} by the correspondence described above.
\end{proof}

\begin{remark}\label{rwkfs}
Using the triangulation of Lemma \ref{lwkfs}, we can realize correspondence 
between the $n$-cubes of \w{W\K} and the \ww{(n+1)}-simplices of the
nerve \w{\N(\K)} simplicially in the barycentric subdivsion $B$ of the
nerve, as follows: 

Consider the triangulated $n$-cube \w{I^{n}\li{\fd}} indexed
by the composable sequence \wref{eqcompseq} as a subcomplex of
\w[.]{\Fs\K} and let \w{B\li{\fd}} denote the barycentric subdivision
of the corresponding \ww{(n+1)}-simplex \w{\sigma^{n+1}\li{\fd}} of
\w[.]{\N(\K)} 

Note that for \w[,]{i\geq 1} the $i$-simplices of 
\w{\sigma^{n+1}\li{\fd}} are labeled by sub-sequences of \w[,]{\fd} with
a single level of parenthesization (indicating where compositions, if
any, have been carried out) - for example, 
\w[.]{(f_{2}f_{3})(f_{4})(f_{5}f_{6})} 
These also label the corresponding vertices of \w{B_{\fd}} (ignoring
those which come from the vertices of \w[),]{\sigma^{n+1}\li{\fd}}
and each $k$-simplex of \w{B\li{\fd}} corresponds to an ascending 
``flag'' of \w{k+1} inclusions of faces of \w[.]{\sigma^{n+1}\li{\fd}}

Now let \w{C\li{\fd}} denote the set of vertices of \w{B\li{\fd}} which are 
labeled by (one-level) parenthesizations of the \emph{full} sequence 
\w{(f_{1},\dotsc,f_{n+1})} (corresponding to the simplices of 
\w{\sigma^{n+1}_{\fd}} which have both \w{a_{0}} and \w{a_{n+1}} as
vertices), and let \w{E\li{\fd}} be the subcomplex of \w{B\li{\fd}}
spanned by \w[.]{C\li{\fd}} 
A $k$-simplex of \w{E\li{\fd}} thus corresponds to a sequence of 
\w{k+1} parenthesizations of \w[,]{(f_{1},\dotsc,f_{n+1})} each
obtained from the next by coalescing a neighbouring pair of
parentheses (since this describes the only face maps of \w{\N(\K)} which
remain inside \w[).]{C\li{\fd}} Therefore, such a $(k+1)$-flag can be
labeled by a single \ww{(k+1)}-level parenthesization
of \w[,]{(f\li{1},\dotsc,f\li{n+1})} just like the $(k-1)$-simplices of 
\w[.]{I^{n}\li{\fd}} Thus \w{I^{n}\li{\fd}} is isomorphic as a simplicial
complex to \w[.]{E\li{\fd}} 
\end{remark}

%
%
\sect{Andr\'{e}-Quillen cohomology of \PAa s}
\label{ccohpa}

Since the category \w{s\PAAlg} of simplicial \PAa s is a
semi-spherical model category (\S \ref{dssmc}), we can use \S
\ref{dcoh} to define the cohomology groups of a \PAa\ $\Lambda$ with
coefficients in a $\Lambda$-module $M$ (see \S \ref{egssmc}).  

\begin{notn}\label{naqcoh}
In such algebraic settings, this is traditionally called 
\emph{Andr\'{e}-Quillen cohomology}, since it can be computed via a
cotangent complex, as in \cite{AndrM,QuiC}. We therefore denote it by 
\w[.]{\HAQ{\ast}{\Lambda}{M}:=\coH{\ast}{\BL}{M}}
\end{notn}

We  would like to compare this with the \ww{\SaO}-cohomology of a
suitable \ww{\SaO}-category (cf.\ \S \ref{ssoc}), for which we need the
following framework: 

\begin{defn}\label{dfl}
Given a set $\A$ of spherical objects in a model category $\C$, we let
\w{\CA} denote the smallest full subcategory  of $\C$ containing $\A$
and closed under weak equivalences and arbitrary coproducts.   

Using \wref[,]{eqsmall} we see that the functor
\w{\piA:\ho\C\to\PAAlg} induces an equivalence of categories  
between the corresponding subcategory \w{\ho\CA} of the homotopy category 
\w{\ho\C} and the category \w{\FA} of \emph{free} \PAa s in \w{\PAAlg}
(namely, those which are isomorphic to \w{\piA B} for \w[).]{B\in\CA}
Moreover, we can extend any \PAa\ \w{\Lambda:\PiA\to\Seta} to a
functor \w{\ho\CA\to\Seta} taking (arbitrary) coproducts
to products.

A \emph{small \ww{\FA}-variant} is a full small subcategory \w{\D} of \w{\FA} 
(or \w[)]{\ho\CA} containing an isomorphic copy of \w[:]{\PiA} in other words, 
\w{\Obj\D} must contain all \emph{finite} coproducts of objects from $\A$, up to
isomorphism.

Given a \PAa\ $\Lambda$ and a small \ww{\FA}-variant $\D$ with
\w[,]{\OO:=\Obj(\D)} we let \w{\Dpp} denote the category with object set
\w[,]{\Op:=\OO\cup\{\bstar\}} where:
\begin{myeq}\label{eqdpp}
\Hom_{\Dpp}(A,B)~=~\begin{cases}
\Hom_{\D}(A,B) & \hsm\text{if}\hsm a,b\in\OO\\
\Hom_{\PAAlg}(A,\Lambda)=\Lambda\lin{A} & \hsm\text{if}\hsm A\in\OO 
\hsm\text{and}\hsm B=\bstar\\
\{\Id_{\bstar},\ast\} & \hsm\text{if}\hsm A=B=\bstar\\
\{\ast\} & \hsm\text{otherwise.}
\end{cases}
\end{myeq}
\noindent That is, all maps out of $\bstar$ are trivial.
Thus we have a full and faithful embedding of $\D$ in \w[,]{\Dpp} 
and $\bstar$ is a weakly terminal object in \w[.]{\Dpp}
We call \w{(\Dpp,\D)} a $\Lambda$-\emph{pair} (in \w[).]{\ho\C} 

Equivalently, if we embed $\D$ in \w{\OpC} (making all maps into
$\bstar$ trivial), we can think of a $\Lambda$-pair (in \w{\ho\C} as
an \ww{\Op}-catgeory under $\D$ (and require only the last three
conditions of \wref[).]{eqdpp} 
\end{defn}
 
\begin{example}\label{egfl}
Let $\D$ be the subcategory of \w{\ho\C} whose objects are of the form
\w{\coprod_{i\in S}~A_{i}} with \w{A_{i}\in\A} \wb{i\in S} and cardinality
\w[.]{\Card(S)\leq\max\{\aleph_{0},\Card(U\Lambda)\}} 
This is a small \ww{\FA}-variant.
We can think of \w{\Dpp} as a subcategory of \w[,]{\PAAlg} by
identifying $\bstar$ with $\Lambda$. 
\end{example}

It turns out that the relative \ww{\SaO}-cohomology of such a pair 
(cf.\ \S \ref{rrelcoh}) has an algebraic interpretation: 

%
%
\begin{thm}\label{tpaso}
Let $\Lambda$ be a \PAa, $M$ a $\Lambda$-module, and \w{(\Dpp,\D)} a
$\Lambda$-pair.  Then for any \w[,]{n\geq 1}
the $n$-th Andr\'{e}-Quillen cohomology group \w{\HAQ{n}{\Lambda}{M}} 
is naturally isomorphic to the $n$-th relative \ww{\SaO}-cohomology group 
\w[.]{\HSO{n}{\Dpp,\D}{M}} 
\end{thm}

\begin{proof}
Let \w{\Vd\to\Lambda} be the canonical free simplicial resolution 
(in the resolution model category on \w{s\PAAlg} of \S \ref{egrmc}(b))
produced by the ``free on underlying'' comonad \w[,]{\F=FU} and let
\w{\Ed} be the analogous free \ww{\SaO}-resolution for 
\w{\Dpp=\D\cup\{\bstar\}} as in \S \ref{egfl}. Thus \w{\pi_{0}\Ed} is 
\w[.]{\Dpp\in\SaOpC}
The relative version \w{\hEd} is obtained from \w{\Ed} by ``excision
of $\D$'' \wh that is, we define the simplicial mapping spaces for
\w{\hEd} by:
$$
\hEd(A,B)~:=~
\begin{cases}
\Ed(A,B) & \text{~if~} B=\bstar\\
\co{\Hom_{\FA}(A,B)} & \text{~if~} B\in \FAp~.
\end{cases}
$$
The twisting map \w{p:\Ed\to\co{\K}} induces an \ww{\SaO}-functor
\w[.]{\rho:\Ed\to\hEd} 

Note that \w{\pi_{0}\Vd\cong U\Lambda} and 
\w[,]{U\Vd\cong\coprod_{\varphi\in\pi_{0}\Vd}~V[\varphi]} where 
\w{V[\varphi]} is the component of \w{\varphi\in\Lambda\lin{A}} for some
\w{A\in\PiA} (\S \ref{dmodels}). Then each \w{V[\varphi]} is isomorphic
to the component of \w{\varphi:A\to\Lambda} in the simplicial mapping space 
\w{\Ed(A,\Lambda)=\hEd(A,\Lambda)} (so in simplicial dimension $n$,
\w{V[\varphi]_{n}} consists of depth $n$ parenthesizations of
composable sequences of morphisms in $\K$, with 
composite $\varphi$).   Because \w{\Vd} is a simplicial \PAa, for any
\w{\theta:A'\to A} in \w{\F'_{\A}} we have a simplicial map 
\begin{myeq}\label{eqaction}
\theta^{\ast}:V[\varphi]\to V[\varphi\circ\theta]
\end{myeq}
\noindent defining an action of \w{\F'_{A}} on the simplicial sets 
\w[.]{V[-]}

Since the catgeory \w{s\PAAlg} is semi-spherical, for each 
$\Lambda$-module $M$ and \w{n\geq 1} we have an 
Eilenberg-Mac~Lane object \w{\EL{M}{n}} in \w[,]{s\PAAlg/\BL} as well 
as an object \w{\EL{D}{n}} in \w[.]{\SaOC/\K} Moreover, we can assume
that both are strict abelian group objects in their respective
categories (see \cite[\S 3.14]{BJTurR}).

Any map of simplicial \PAa s \w{f:\Vd\to\EL{M}{n}} (over \w[)]{\BL}
defines an \ww{\SaO}-map \w[,]{\hat{f}:\hEd\to\EL{D}{n}} which is
defined on the simplicial mapping spaces \w{\Ed(A,\Lambda)} via
the above identification with the components \w{V[\varphi]} of
\w[.]{\Vd} These fit together to define an \ww{\SaO}-map, because of
the action \wref[.]{eqaction}

Precomposing this with \w{\rho:\Ed\to\hEd} yields an element in the
relative \ww{\SaO}-cohomology group \w[.]{\HSO{n}{\Dpp,\D}{M}} Similarly
for the converse direction.  
\end{proof}

%
%
\sect{Diagrams and \ww{\SO}-categories}
\label{cdoc}

We now explain the approach of Dwyer, Kan, and Smith (see
\cite{DKanS,DKanO,DKSmiO,DKSmiR}) to realizing a homotopy-commutative
diagram \w[,]{X:\K\to\ho\TT} based on the concepts introduced in
Section \ref{csimpmc}. 

\begin{defn}\label{dhhc}
A \emph{diagram up to homotopy} in a simplicial model category $\C$ 
is a functor \w{X:\K\to\ho\C} from some small indexing category 
$\K$. By definition, one can choose a functor \w{X_{0}:\sk{0}F(\K)\to\C} 
lifting $X$ (sometimes called a \emph{$0$-realization} of $X$). 
An extension of any such a \w{X_{0}} to a simplicial functor 
\w{X_{\infty}:F(\K)\to\C} makes $X$ \emph{$\infty$-homotopy commutative}.  
\end{defn}

A classical result of Boardman and Vogt (compare \cite[Corollary 2.5]{DKSmiR}) 
says:

%
%
\begin{thm}[\protect{\cite[Cor.~4.21 \& Thm.~4.49]{BVogHI}}]\label{ftwo}
A diagram \w{X:\K\to\ho\TT} can be rectified (i.e., lifted to 
\w[)]{\hX:\K\to\TT} if and only if $X$ can be made $\infty$-homotopy
commutative. 
\end{thm}

\begin{notn}\label{nsec}
When we want to emphasize that we are thinking of a simplicial model category $\C$ 
just as a simplicially enriched category, we denote it by \w[.]{\sC}
\end{notn}

\begin{remark}\label{rcx}
Theorem \ref{ftwo} implies that the rectification of a homotopy
commutative diagram \w{X:\K\to\ho\C} can be described in purely in
terms of the simplicially enriched category \w{\sC} \wh in fact, we
can restrict to an \ww{\SaO}-category \w[,]{\CX} the sub-simplicially
enriched category of \w{\sC} with function complex
\w{\map_{\CX}(u,v):=\map_{\sC}(Xu,Xv)} for each \w[.]{u,v\in\OO:=\Obj\K} 

Note that a choice of a $0$-realization \w{X_{0}:\Gamma\to\Ta}
is equivalent to choosing basepoints in each \w[,]{\CX(u,v)}
though of course this cannot be done coherently unless $X$ is rectifiable.
\end{remark}

\begin{mysubsection}[\label{sdks}]{The obstruction theory}
Given \w{X:\K\to\ho\C} as above, the (possibly empty) moduli space 
\w{\hc{}X} of all rectifications of $X$ is homotopy equivalent to the
space \w{\hc{\infty}X:=\map_{\OC}(F(\K),\CX)} of all functors making $X$
$\infty$-homotopy commutative, which in turn is the (homotopy)
inverse limit of the tower:
$$
\hc{\infty}X \to \dotsc \to \hc{n}X\to \hc{n-1}X\dotsc \to \hc{1}X~,
$$
\noindent where \w[.]{\hc{n}\A:=\map_{\OC}(F(\K),\Po{n-1}\CX)} Therefore,
the realization problem can be solved if one can successively lift 
\w{\hX_{1}\in\hc{1}X} through the tower. 

The components of \w{\hc{\infty}X} are not in general determined by
those of  the spaces \w{\hc{n}X} (cf.\ \cite[3.4]{DKSmiR}).
Because each \w{\hc{n}X} is a mapping space, we can use successive
liftings \w{\hX_{n}\in \hc{n}X} to pull back the \ww{(n-1)}-st $k$-invariant
for \w{\CX} to a map \w[,]{h_{n}:\Fs\K\to K^{G}(\pi_{n}\CX,n+1)} and
Dwyer, Kan, and Smith show:
\end{mysubsection}

%
%
\begin{prop}[\protect{\cite[Proposition 3.6]{DKSmiR}}]\label{fsix}
The map \w{\hX_{n}} lifts to \w{\hX_{n+1}\in\hc{n+1}X} if and only
if \w{[h_{n}]\in\HSO{n+1}{\K}{\pi_{n}\CX}} vanishes.
\end{prop}

\begin{mysubsection}[\label{srelvo}]{A relative version}
There is also a relative version of this obstruction theory, in which, 
given \w{X:\K\to\ho\C} as above, we assume that we have a subcategory $\LL$ 
of $\K$ equipped with a lift \w{\hY:\LL\to\C} of \w[.]{X\rest{\LL}}
This defined a map from the pushout 
$$
\xymatrix@R=25pt{
\ar @{}[drr] |>>>>>{\framebox{\scriptsize{PO}}}
\Fs\LL \ar[d]_{p} \ar[rr]^{i} & & \Fs\K \ar[d]\\
\LL \ar[rr] & & \Fs(\K,\LL)~,
}
$$
\noindent (compare \wref[)]{eqpo} into \w[,]{\Po{0}\CX} which lifts 
(non-canonically) to \w{\Po{1}\CX} because $X$ is homotopy-commutative.

Again, we can use each successive liftings \w{\hX_{n}\in \hc{n}X} to pull back 
the \ww{(n-1)}-st $k$-invariant for \w{\CX} to a map  
\w[,]{h_{n}:\Fs(\K,\LL)\to K^{\K}(\pi_{n}\CX,n+1)} representing an 
\ww{\SaO}-cohomology class \w[,]{[h_{n}\in\HSO{n+1}{(\K,\LL)/\K}{\pi_{n}\CX}}
and the relative version of Proposition \ref{fsix} clearly holds.
\end{mysubsection}

%
%
\sect{The first obstruction}
\label{cfobst}

Given a natural system $D$ on a category $\K$, one can always construct a 
trivial linear track category with $D$ as its (abelian) fundamental groupoid.
Moreover, by Proposition \ref{tlintrext}, the linear track extensions of $\K$
by $D$ are classified up to weak equivalence by \w[.]{\HBW{3}{\K}{D}}
When \w{\K=\PiA} for some set $\A$ of spherical objects in a model category 
$\C$, the cohomology class determining the extension is represented by the 
explicit cocycle of \S \ref{rcocycle}.  We now show how this is reflected 
in \w[.]{\SaOC} For this purpose, we need an \ww{\SaO}-version of 
Definition \ref{dfl}:

\begin{defn}\label{dcap}
Let $\C$ be a simplicial model category with spherical objects $\A$.
A \emph{small \ww{\sCA}-variant} is a full (necessarily small) fibrant 
sub-simplicial category \w{\hCA} of \w{\sCA} (\S \ref{nsec}), such that 
\w{\pi_{0}\hCA} is a small \ww{\FA}-variant (\S \ref{dfl}). This just means 
that \w{\OO:=\Obj\hCA} contains all finite coproducts of objects of $\A$, 
up to weak equivalence. 

We assume that all objects in \w{\CA'} are cofibrant, and for simplicity we
also assume that $\OO$ contains a \emph{canonical} copy of \w[.]{\Obj\PiA}
\end{defn}

\begin{example}\label{egcap}
A \emph{minimal} small \ww{\sCA}-variant is any skeletal subcategory $\fX$ of 
\w{\sCA} with \w{\pi_{0}\fX=\PiA} (\S \ref{dmodels}). In particular, we denote
by \w{\CAm} the \emph{canonical} minimal small \ww{\sCA}-variant, whose objects
consist of a (functorial) fibrant and cofibrant replacement for each non-isomorphic 
finite coproduct of objects from $\A$. 

More generally, if $\D$ is any small \ww{\FA}-variant, choose any 
embedding \w{i:\D\hra\ho\C} for which \w{i(a)} is fibrant and cofibrant for 
each \w[.]{a\in\OO:=\Obj\D} We then obtain a fibrant small \ww{\sCA}-variant 
\w{\hCA} by setting \w[.]{\map_{\hCA}(a,b):=\map_{\CA}(i(a),i(b))}
\end{example}

\begin{mysubsection}[\label{szerok}]{The $0$-th $k$-invariant}
In general, it makes no sense to speak of the $0$-th $k$-invariant of 
an \ww{\SaO}-category $\X$, since \w{\hpi\X} is not an abelian group 
object over \w{\K:=\pi_{0}\X} \wh even though we do have a pullback square
of the form \wref{eqkinv} for \w[,]{n=0} too. However, \w{k_{0}} is a
well-defined cohomology class in the following specialized situation:
\end{mysubsection}

\begin{assume}\label{alinext}
Let $\A$ be a collection of spherical objects in a simplicial model
category \w[,]{\sC} let \w{\hCA} be a small \ww{\sCA}-variant \wh so that 
\w{\D:=\pi_{0}\hCA} is a small $\FA$-variant),
the track category $\E$ of \w{\hCA} is linear (\S \ref{egtrackcat}), and 
\w{\Omega\D} is a natural system on \w{\D:=\pi_{0}\hCA} (cf.\ \S \ref{emodule} 
and \S \ref{egnatsys}). We let \w[.]{\OO:=\Obj\hCA}
\end{assume}

With these assumptions we find:

%
%
\begin{thm}\label{tzero}
The $0$-th $k$-invariant for \w{\hCA} corresponds to the 
cohomology class \w{\chi_{\E}} classifying the linear track extension 
\w{\Omega\D\to\E\to\D}  (cf.\ \S \ref{rcocycle}) under the natural 
isomorphism of Theorem \ref{tbwso}. 
\end{thm}

\begin{proof}
Set \w[,]{\X:=\hCA} and consider the following square of the form 
\wref{eqkinv} in \w[:]{\SaO}
$$
\xymatrix@R=25pt{
\ar @{}[drr] |>>>>>{\framebox{\scriptsize{PO}}}\Po{1}\X \ar[d]_{p}\ar[rr]^{i} 
& & \Y \ar[d]^{q} & \Fs\D=B\D \ar[l]_{\xi}^{\simeq} \\
\Po{0}\X \ar[rr]_{j} & & \Z~, & 
}
$$
\noindent in which the homotopy pushout \w{\Z} satisfies 
\w{\Po{2}\Z\simeq E_{\D}(\Omega\D,2)\simeq E_{\D}(\pi_{1}\X,2)}
by \cite[Proposition 6.4]{BDGoeR}, and thus if \w{r\q{2}:\Z\to\Po{2}\Z}
is the structure map of \S \ref{dssmc}(c), the $0$-th $k$-invariant for
$\X$ is \w{k_{0}:=r\q{2}\circ q\sim r\q{2}\circ j} by construction.
We use Kan's orginal model for the Postnikov system, so that
\w{(\Po{k}X)_{n}} consists of \ww{\sim_{k}}-equivalence classes of
$n$-simplices in \w[,]{\X} where 
\w[,]{\sigma\sim_{k}\tau\ \EQUIV\ \sk{k}\sigma=\sk{k}\tau}  
(cf.\ \cite[VI, \S 2]{GJarS}). We assume that \w{\X} is fibrant (so
each mapping space \w{X(u,v)} is a Kan complex).

Factor \w{p=p\q{1}:\Po{1}\X\to\Po{0}\X} as a cofibration 
\w{i:\Po{1}\X\to\Y} followed by a weak equivalence, so that the
pushout above is a homotopy pushout, as required. 
Thus \w{Y_{i}=X_{i}} for \w[,]{i\leq 1} while
\w[,]{Y_{2}=(X_{2}/\sim)\,\amalg\,\bar{Y}} where \w{\bar{Y}} has a 
``fill-in'' $2$-simplex \w{T=T_{(\sigma_{0},\sigma_{1},\sigma_{2})}}
for every triple of $1$-simplices \w{(\sigma_{0},\sigma_{1},\sigma_{2})} 
in \w[,]{X_{1}} with matching faces, having \w[.]{d_{i}T=\sigma_{i}} 
The pushout \w{\Z} thus consists of the reduction via \w{\sim_{0}} of
the copy of \w{\X} in \w[,]{\Y} with $\bar{Y}$ unaffected. The
$2$-simplices \w{K_{(\sigma,0,0)}} for non-null homotopic $\sigma$
represent \w{\hpi\X} in \w[.]{\Po{2}\Z\simeq E_{\D}(\hpi\X,2)}
We shall not need the description of \w{\Y} or \w{\Z} in higher
dimensions. 

Let \w{\Fs\D} be the cofibrant replacement for \w{\Po{0}\X}
constructed as in \S \ref{ssoc}. The weak equivalence \w{\xi:\Fs\D\to\Y}
is then defined as follows: 

Every $0$-simplex \w{(\phi)\in\Fs\D} corresponds to a homotopy
class \w[,]{\phi\in[Xu,Xv]_{\ho\C}} and \w{\xi(\phi)} is a choice of a 
representative \w{s(\phi)} in \w[.]{(\Po{0}\X)_{0}=X(u,v)_{0}} For a
(non-composite) $1$-simplex \w{\sigma=((\phi_{k})\dotsc(\phi_{1}))} in
\w[,]{(FU)^{2}\D} \w{\xi(\sigma)} is a choice of a homotopy 
\w{H_{(\phi_{k},\dotsc,\phi_{1})}} between 
\w{s(\phi_{k})\cdot \dotsc\cdots(\phi_{1})} and 
\w[,]{s(\phi_{k}\cdot\dotsc\cdot\phi_{1})} which exists since 
\w[.]{\D=\ho\E} Finally, the faces of any $2$-simplex
\w{\tau\in(FU)^{3}\D} form a triple of matching $1$-simplices, so their
image under $\xi$ has a canonical fill-in \w[,]{T\in\bar{Y}} and we
set \w[.]{\xi(\tau)=T} 

Now either of the two maps from \w{\Po{0}\X} to \w{\Po{2}\Z} represents
\w[;]{k_{0}} using the cofibrant model \w{\Fs\D} for the source, it
is enough to identify the map on $2$-simplices \wh or, using the
identification of simplicial and cubical cohomology mentioned in the
proof of Theorem \ref{tbwso}, on the (triangulated) square 
\w{I^{2}\li{\phi_{3}\circ\phi_{2}\circ\phi_{1}}} as in \wref[.]{eqsquare}
By the descriptions of $\xi$ and $\Y$ above, this maps to:

\begin{center}
%
%
\begin{picture}(300,180)(40,0)
%
%
\put(-8,152){$s(\phi_{3})\cdot s(\phi_{2})\cdot s(\phi_{1})$}
\put(95,155){\circle*{5}}
\put(293,155){\vector(-1,0){190}}
\put(160,161){$H_{(\phi_{3},\phi_{2})\cdot s(\phi_{1})}$}
\put(295,155){\circle*{5}}
\put(302,152){$s(\phi_{3}\phi_{2})\cdot s(\phi_{1})$}
%
%
\put(95,20){\vector(0,1){130}}
\put(18,82){$s(\phi_{3})\cdot H_{(\phi_{2},\phi_{1})}$}
\put(95,15){\circle*{5}}
\put(13,12){$s(\phi_{3})\cdot s(\phi_{2}\phi_{1})$}
%
%
\put(290,15){\vector(-1,0){190}}
\put(170,0){$H_{(\phi_{3},\phi_{2}\phi_{1})}$}
\put(295,15){\circle*{5}}
\put(302,12){$s(\phi_{3}\phi_{2}\phi_{1})$}
%
%
\put(295,20){\vector(0,1){130}}
\put(300,82){$H_{(\phi_{3}\phi_{2},\phi_{1})}$}
%
%
\put(183,96){\vector(-3,2){82}}
\put(170,80){$H_{(\phi_{3},\phi_{2},\phi_{1})}$}
\put(210,74){\line(3,-2){79}}
%
%
\put(180,123){$\xi(((\phi_{3})(\phi_{2}))((\phi_{1})))$}
\put(118,40){$\xi(((\phi_{3}))(\phi_{2})(\phi_{1})))$}
\end{picture}
\end{center}

\noindent which is just the cocycle of \wref[,]{eqbauescocyc} under the
isomorphism of Theorem \ref{tbwso}.
\end{proof}

\begin{cor}\label{czero}
Under the assumptions of \S \ref{alinext}, the equivalence classes of
linear track extensions \w{\Omega\D\to\E\to\D} are in one-to-one
correspondence with one-stage Postnikov systems of \ww{\SaO}-categories
$\Y$ (that is, those satisfying \w[)]{\Y\simeq\Po{1}\Y} such that 
\w[,]{\pi_{0}\Y\cong\D} and \w{\pi_{1}\Y\cong\Omega\D} as
\ww{\iO{\RM{\K}}}-categories.
\end{cor}

\begin{mysubsection}[\label{srelv}]{A relative version}
Now assume that \w{\hCA\in\SaOC} as in \S \ref{alinext} extends to an 
\ww{\SaOp}-subcategory $\X$ of \w[,]{\sC} obtained by adding a single new
object \w[.]{Y\in\C} Thus \w[,]{\Op:=\OO\cup\{Y\}} \w[,]{\X\rest{\OO}=\hCA}
and we omit all non-trivial maps out of $Y$, so that
\w{\map_{\X}(Y,Y)=\co{\{\Id_{Y},\ast\}}} and
\w{\map_{\fX}(Y,B)=\co{\{\ast\}}} for all \w{B\in\OO} (see \S
\ref{dmapa} below).

In this case we can extend the track category $\E$ of \w{\hCA} to a
track category \w{\Ep} for $\X$, which is still linear (since all 
non-trivial maps are out of homotopy cogroup objects). If 
\w[,]{\Dpp:=\pi_{0}\X} then \w{(\Dpp,\D)} is a $\Lambda$-pair, for
\w{\Lambda:=\piA Y} (\S \ref{dfl}), and   \w{\Omega\Dpp} is a natural
system on \w[.]{\Dpp}  Therefore, Theorem \ref{tzero} applies in this
situation, too: that is, the $0$-th $k$-invariant for \w{\X}
corresponds to the cohomology class classifying the linear track
extension \w[.]{\Omega\Dpp\to\Ep\to\Dpp}  

Note that the inclusion of categories \w[,]{\hCA\hra\X} and the
corresponding inclusion of objects sets \w[,]{\OO\hra\Op} induces 
natural transformations in Baues-Wirsching and \ww{\SaO}-cohomology 
fitting into long exact sequences with the relative versions, with all
vertical maps being isomorphisms by Theorem \ref{tbwso}:
\mydiagram[\label{eqles}]{
\dotsc\HBW{n}{\Dpp}{\Omega\Dpp} \ar[r]^>>>>>{i^{\ast}} \ar[d]^{\cong} &
\HBW{n}{\D}{\Omega\D} \ar[r]^>>>>>{\delta^{n}} \ar[d]^{\cong} &
\HBW{n+1}{\Dpp,\D}{\Omega\Dpp}\dotsc \ar[d]^{\cong} \\
\dotsc\HSO{n-1}{\X}{\Omega\Dpp} \ar[r]^>>>>>{i^{\ast}} & 
\HSO{n-1}{\hCA}{\Omega\D}\ar[r]^>>>>{\delta^{n-1}} & 
\HSO{n}{\X,\hCA}{\Omega\Dpp} \dotsc
}
\end{mysubsection}

%
%
\begin{lemma}\label{pfobst}
The class \w{\delta^{3}(\chi_{\E})} in \w{\HSO{4}{\Dpp,\D}{\Omega\Dpp}}
is the obstruction to realizing $\Lambda$ by a track category
\w{\Ep} inside that of $\C$. 
\end{lemma}

\begin{proof}
The class \w{\delta^{3}(\chi_{\E})} vanishes if and only if
\w{\chi_{\E}} is in the image of \w{i^{\ast}} in the top row of
\wref{eqles} \wh that is, if and only if \w{\E} extends to a linear
track category \w{\Ep} realizing $\Lambda$.
\end{proof}

From the ladder of isomorphisms \wref{eqles} we deduce:

\begin{cor}\label{cofobst}
The class \w{\delta^{2}(\chi_{\tD})} maps under the isomorphism of
Theorem \ref{tbwso} to the relative $k$-invariant
\w{\delta^{2}(k_{0})} in \w[,]{\HSO{3}{\X,\hCA}{\Omega\Dpp}}
which is the obstruction to realizing $\Lambda$ as a one-stage
Postnikov system in the \ww{\SaO}-category for $\C$. 
\end{cor}

%
%
\sect{Realizing \PAa s}
\label{crpa}

The approach of \cite{DKStB,DKStB,BDGoeR} to realizing \Pal s 
can be generalized somewhat (see \cite{BJTurR}), but it still does not
apply to arbitrary resolution model categories (for example, it does
not even apply to topological spaces, if $\A$ consists of mod-$p$
Moore spaces \wh see \cite[\S 4.6]{BlaM}). We therefore restrict to
the following setting: 

\begin{defn}\label{demc}
If $\C$ is a semi-spherical resolution model category equipped with a
set of spherical objects $\A$, the resolution model category \w{s\C}
(\S \ref{drmc}) is called a \emph{strict \ww{E^{2}}-model category} if
the inclusion \w{\co{-}:\C\to s\C} has a left adjoint
\w[,]{R:s\C\to\C} called the \emph{realization} functor for \w[,]{s\C}
such that for all \w[,]{A_{0}\in\A} the natural map induced by the unit 
\w[:]{\var_{\Xd}:\Xd\to\co{R\Xd}} 
\begin{myeq}[\label{eqrealmap}]
\var_{\ast}:\|\map_{\C}(A_{0},\Xd)\|\to \map_{\C}(A_{0},R(\Xd))
\hsm\text{is a weak equivalence}
\end{myeq}
\noindent as long as \w{\Xd\in s\C} is cofibrant in the resolution
model category structure on \w{s\C} determined by
\w[.]{\A_{0}:=\{\Sigma^{k}A_{0}\}_{k=0}^{\infty}} Here \w{\|\Qd\|} is
the diagonal of a bisimplicial set \w[.]{\Qd\in s\Ss} 
\end{defn}

\begin{example}\label{egemc}
The main example we have in mind is \w{\C=\Ta} with
\w[,]{\A=\{\bS{k}\}_{k=1}^{\infty}} and $R$ the usual geometric
realization. In this case the cofibrancy condition on \w{\Xd} implies
that each \w{X_{n}} is \ww{(k-1)}-connected, when \w[,]{A_{0}=\bS{k}}
so \wref{eqrealmap} holds by \cite[Theorem 12.3]{MayG} (see also
\cite{AndF}). 

In \cite[Theorems 3.15-3.19]{BJTurR}, it was shown that all the
examples of \S \ref{egrmc} are \ww{E^{2}}-model categories, which
satisfy a somewhat weaker set of axioms (see \cite[Definition
  3.12]{BJTurR}).  However, there are a number of additional examples
satisfying these stricter conditions \wh \w{\Ta} can be replaced by
\w{\Sr} or $\G$, or various categories of spectra, or DG-categories;
or we can take diagrams in these categories. We can also use localized
or truncated spheres. In order to cover all these cases we have
therefore  stated the conditions needed in axiomatic form. This also
permits them to be dualized more readily (\S \ref{rdpa}).
\end{example}

In this context the obstruction theory of \cite{BDGoeR} can be stated
using the following

\begin{defn}\label{dquap}
A \emph{quasi-Postnikov tower} for a \PAa\ $\Lambda$ is a tower of fibrations:
$$
\dotsb\xra{p\q{n+1}}\Xn{n+1}\xra{p\q{n}}\Xn{n}\xra{p\q{n-1}}\dotsb\xra{p\q{0}}
\Xn{0}\simeq\BL
$$
\noindent in \w{s\C/\BL} such that 
\w{\piA\Xn{n}\simeq\EL{\Omega^{n+1}\Lambda}{n+2}} (as for the
usual Postnikov system of a realization of \w{\BL} in \w{s\C} \wh see
\cite[\S 5.8]{BJTurR}). The object \w{\Xn{n}\in s\C} will be called an 
\emph{$n$-th quasi-Postnikov section} for $\Lambda$. 
\end{defn}

The following is shown in \cite[\S 9]{BDGoeR} and \cite[Theorems 5.6-5.7]{BJTurR}:
%
%
\begin{thm}\label{treal}
Let $\C$ be an \ww{E^{2}}-model category with a set of spherical objects $\A$.
A \PAa\ $\Lambda$ is realizable if and only it has a quasi-Postnikov tower 
in \w[.]{s\C/\BL} Moreover, if such a tower exists in degrees \w[,]{\leq n-1}
then: 
\begin{enumerate}
\renewcommand{\labelenumi}{(\alph{enumi})~}
\item Up to homotopy, there is a unique \w{\Xn{n}\in s\C} with 
\w[,]{\Po{n-1}\Xn{n}=\Xn{n-1}}
\begin{myeq}[\label{eqpnpsys}]
\pinat{k}{\Xn{n}}\cong\begin{cases}
		\Omega^{k}\Lambda &~~\text{for \ } 0\leq k\leq n,\\
		 0     &~~\text{otherwise}~\end{cases}
\end{myeq}
\noindent (see \S \ref{emodule}). 
\item  This \w{\Xn{n}} is an $n$-th quasi-Postnikov section for
  $\Lambda$ if and only if the \ww{(n+2)}-nd $k$-invariant for
  \w{\piA\Xn{n}} vanishes in
  \w[.]{\HAQ{n+3}{\Lambda}{\Omega^{n+1}\Lambda}}  
\item In that case, the different choices for the map
  \w{p\q{n}:\Xn{n+1}\to\Xn{n}} are in one-to-one correspondence with
  elements of \w[.]{\HAQ{n+2}{\Lambda}{\Omega^{n+1}\Lambda}} 
\end{enumerate}
\end{thm}

Note that from the spiral exact sequence \wref{eqspiral} 
we can deduce from \wref{eqpnpsys} that
\begin{myeq}[\label{eqethg}]
\pi_{k}\piA\Xn{n}\cong\begin{cases}
		\Lambda &~~\text{for \ } k=0\\
		\Omega^{n+1}\Lambda &~~\text{for \ } k=n+2,\\
		 0     &~~\text{otherwise}~\end{cases}
\end{myeq}
\noindent The vanishing of the \ww{(n+2)}-nd $k$-invariant for \w{\piA\Xn{n}}
is equivalent to the latter being an Eilenberg-Mac~Lane object 
\w{\EL{\Omega^{n+1}\Lambda}{n+2}} in \w[.]{s\PAAlg}

%
%
\sect{Mapping algebras}
\label{cma}

In order to compare the approach of Sections \ref{cdoc} and \ref{crpa}, 
we need to recast the problem of realizing \w{\Lambda\in\PAAlg} as one of 
rectifying a suitable homotopy-commutative diagram \wh or more precisely, of
lifting a diagram through the Postnikov system of an \ww{\SaO}-category.

The obvious first choice is to consider a diagram \w{X:\K\to\ho\C} for
\w{\K:=\FAp} (\S \ref{dfl}). Unfortunately, there are two problems with this:    

\begin{enumerate}
\renewcommand{\labelenumi}{(\alph{enumi})}
\item We do not actually have such a diagram $X$ to begin with, since the
putative value of \w{X(\bstar)\in\ho\C} is precisely the realization of
the \PAa\ $\Lambda$ in $\C$ that we are looking for. 
\item Moreover, we do not expect a rectification \w{\hX:\FAp\to\C} to
  exist (unless the model category $\C$ is ``formal''),
since commuting diagrams in \w{\ho\C} do not generally lift to $\C$ 
\end{enumerate}

In order to solve the second problem, we introduce the following concept:

\begin{defn}\label{dmapa}
Let \w{\hCA} be a small \ww{\sCA}-variant (\S \ref{dcap}) with object set $\OO$, 
and let \w[.]{\Op:=\OO\cup\{\bstar\}}
An \emph{\Ama\ (based on \w{\hCA})} is an \ww{\SaOp}-category 
\w{\fX} with mapping spaces as follows (compare \wref[):]{eqdpp}
\begin{myeq}\label{eqhca}
\map_{\fX}(B,C)~=~\begin{cases}
\map_{\sCA}(B,C) & \hsm\text{if}\hsm B,C\in\OO\\
\co{\{\ast,\Id_{\bstar}\}} & \hsm\text{if}\hsm B=C=\bstar\\
\co{\{\ast\}} & \hsm\text{otherwise}
\end{cases}
\end{myeq}

The category of all \Ama s based on \w{\hCA} will be denoted by 
\w{\MAC} (or simply \w[,]{\MA} when \w{\hCA} is understood from the context). 
Elements in \w{\MA} will be written $\fX$, $\fY$, etc, and we denote 
\w{\map_{\fX}(B,\bstar)} by \w{\fX\lin{B}} for all \w[.]{B\in\OO}  
If we embed \w{\SaOC} in \w{\SaOpC} by making
\w{\map(B,\bstar)=\{\ast\}} for all \w{B\in\Op} (as in \S \ref{dfl}),  
then we can think of an \Ama\ based on \w{\hCA} as an
\ww{\SaOp}-category under \w[,]{\hCA} subject to last two conditions
of \wref[.]{eqhca}  Thus \w{\MA} inherits a simplicial model category
structure from \w[.]{\SaOpC}  

Note that if we set \w[,]{\Dpp:=\pi_{0}\fX} we obtain a $\Lambda$-pair
\w{(\Dpp,\D)} for \w[,]{\D:=\pi_{0}\hCA} where the \PAa\ $\Lambda$ is
defined by \w{\Lambda\lin{A}:=\pi_{0}\fX\lin{A}} for all \w[.]{A\in\A}
Thus we can think of an \Ama\ as an enriched version of a \PAa. 
\end{defn}

\begin{example}\label{egmapa}
Given a small \ww{\sCA}-variant \w[,]{\hCA\subseteq\sCA} the motivating 
example of a \Ama\ $\fX$ based on \w{\hCA} is obtained by choosing any 
\w[,]{X\in\C} and setting \w[.]{\map_{\fX}(A,\bstar):=\map_{\sC}(A,X)}
We denote this \Ama\ by \w{\fMA^{\hCA}X} (or simply \w[,]{\fMA X} when
\w{\hCA} is understood from the context). Clearly 
\w[.]{\pi_{0}(\fMA X)\cong\piA X} We say that an \Ama\ $\fY$ is
\emph{realizable} (by \w[)]{X\in\C} if \w[.]{\fY\cong\fMA X} Since any
\w{Y\in\C} is fibrant, \w{\fMA X} is always fibrant.
\end{example}

\begin{remark}\label{rloops}
Recall that the \emph{path object} \w{PK\in\Sa} for a fibrant pointed
simplicial set $K$ has 
\w[,]{(PK)_{n}:=\{x\in K_{n+1}~:\ d_{1}\dotsc d_{n+1} x=\ast\}} 
with re-indexed face and degeneracy maps, and the universal fibration
\w{p:PK\to K} is induced by \w{d_{0}} (cf.\ \cite[\S 2.9]{CurtS}). 
We denote the \emph{path fibration} functor 
\w{K\mapsto(PK\xra{p}K)} by \w[,]{\rho:\Sa\to\Sa^{T}} where
\w{\Sa^{T}} is the category of diagrams in \w{\Sa} indexed by
\w[.]{T=(0\to 1)} Because $\rho$ commutes with products, it extends to
a functor \w[.]{\rho:\MA\to\MA^{T}} 

Note that:
\begin{myeq}[\label{eqconepath}]
\rho\map_{\sC}(A,Y)\hsp \text{is induced by the inclusion}\hsm i:A\hra CA
\end{myeq}

If we define the suspension \w{\Sigma X} in $\C$ as the cofiber of
\w[,]{i:X\hra CX} where \w{CX} is the reduced cone, then for any
fibrant \Ama\ $\fX$ and every \w{B\in\hCA} we have a natural map
$\zeta$ to the pullback (in \w[),]{\Sa} as indicated: 
\mydiagram[\label{eqloops}]{
\fX\lin{\Sigma B} \ar@/^2pc/[drrr]^{i^{\ast}}
\ar@{.>}[dr]^{\zeta} \ar@/_2pc/[ddr] &&&\\
& \Omega\fX\lin{B} \ar[rr] \ar[d] 
\ar @{} [drr] |<<<{\framebox{\scriptsize{PB}}} && P\fX\lin{B} \ar[d]^{p} \\
& \ast \ar[rr] && \fX\lin{B} 
}

Similarly, if \w{B=\coprod_{i\in I} B_{i}} for \w[,]{B_{i}\in\hCA}
we have a natural map
\begin{myeq}[\label{eqcoprodma}]
\fX\lin{B}~\xra{\theta}~\prod_{i\in I}\fX\lin{B_{i}}
\end{myeq}
\end{remark}

\begin{defn}\label{drealistic}
An \Ama\ $\fX$ based on \w{\hCA} will be called \emph{realistic} if whenever
there are weak equivalences 
\begin{myeq}[\label{eqrealistic}]
A'\simeq\Sigma A  \hsp\text{and}\hsp B\simeq\coprod_{i\in I} B_{i}~,
\end{myeq}
\noindent in \w[,]{\hCA} the maps $\zeta$ in \wref{eqloops} and $\theta$ 
in \wref{eqcoprodma} are weak equivalences, 
\end{defn}

\begin{lemma}\label{lrealist}
Any realizable \Ama\ is realistic.
\end{lemma}

\begin{proof}
This holds since both $\zeta$ in \wref{eqloops} and $\theta$ in 
\wref{eqcoprodma} map into homotopy limits.
\end{proof}

Note that if \w{\fX:=\fMA Y} and one of the maps in \wref{eqrealistic} 
is an isomorphism, so is the corresponding map $\zeta$ or $\theta$.

\begin{lemma}\label{lmapeq}
Any map \w{f:X\to X'} in $\C$ induces a morphism of \Ama s 
\w[,]{f_{\ast}:\fMA Y\to\fMA Y'} and $f$ is an $\A$-equivalence (\S
\ref{dmodels}) if and only if \w{f_{\ast}:\fMA Y\lin{A}\to\fMA Y'\lin{A}} 
is a weak equivalence in \w{\Sa} for each
\w[.]{A\in\A}\hfill $\Box$ 
\end{lemma}

\begin{defn}\label{dfreema}
A \emph{free} \Ama\ based on \w{\hCA} is one of  the form \w{\fMA B} for 
\w[.]{B\in\hCA} 
\end{defn}

\begin{lemma}\label{lfreema}
If $\fY$ is an \Ama\ based on \w{\hCA} and \w[,]{B\in\hCA} there is a
natural isomorphism \w[.]{\map_{\MA}(\fMA B,\fY)\cong\fY\lin{B}}
\end{lemma}

\begin{proof}
This follows from the enriched Yoneda Lemma (cf.\ \cite{DKellE}).
\end{proof}

\begin{defn}\label{dpostma}
If $\fX$ is a \Ama\ based on \w[,]{\hCA} for any \w{n\geq 0} we obtain 
its $n$-\emph{th Postnikov section} \w{\Po{n}\fX} by setting 
\w{(\Po{n}\fX)\lin{B}:=\Po{n}(\fX\lin{B})} for any \w[.]{B\in\OO:=\Obj\hCA}
This is well-defined, since when we compose the composition map 
\w{\gamma:\map_{\fX}(B,A)\times\map_{\fX}(A,\bstar)\to 
\map_{\fX}(B,\bstar)=\fX\lin{B}}
of the simplicial enrichment with Postnikov fibration 
\w[,]{p:\fX\lin{B}\to\Po{n}(\fX\lin{B})} the result factors
as:
\begin{equation*}
\begin{split}
\map_{\hCA}&(B,A)\times\map_{\fX}(B,\bstar)~\to~
\Po{n}\map_{\hCA}(B,A)\times (\Po{n}\fX)\lin{A}\\
&=~\Po{n}\left(\map_{\hCA}(B,A)\times\fX\lin{A}\right)~\xra{\Po{n}\gamma}~
(\Po{n}\fX)\lin{B}~.
\end{split}
\end{equation*}

A map of \Ama s \w{\Phi:\fX\to\fY} is called an \emph{$n$-equivalence} 
if it induces a weak equivalence of $n$-th Postnikov sections. 
A map \w{f:X\to Y} in $\C$ is an \emph{$n$-stage $\A$-equivalence} if
\w{\fMA f:\fMA X\to \fMA Y} is an \emph{$n$-equivalence} of \Ama s.
\end{defn}

\begin{example}\label{egpostma}
For \w[,]{n=0} we can replace \w{\Po{0}\fX} by \w[,]{\pi_{0}\fX} using the fact 
that the composition in any simplicially enriched category $\X$ factors through 
its homotopy category \w[.]{\pi_{0}\X} In particular, this shows that
if \w{\hCA} is a small \ww{\sCA}-variant and \w[,]{\D:=\pi_{0}\hCA}
then any $\Lambda$-pair \w{(\Dpp,\D)} can be enriched by an \Ama\
\w{\fX_{\Lambda}} based on \w{\hCA} with \w[.]{\pi_{0}\fX_{\Lambda}\cong\Lambda}
\end{example}

\begin{remark}\label{rpostma}
The tower \w{(\Po{n}\fMA X)_{n=0}^{\infty}} may be the best
approximation to an $\A$-Postnikov tower available, since the category
$\C$ itself may not have such towers \wh e.g., when \w{\C=\Ta} and
$\A$ consist of mod-$p$ Moore spaces (see \cite[\S 3.10]{BlaC}).
\end{remark}

%
%
\sect{The Stover category}
\label{csc}

We now specialize to a specific small \ww{\sCA}-variant, which defines 
a kind of \Ama s with various useful properties:

\begin{defn}\label{dstovcat}
Let $\C$ be an \ww{E^{2}}-model category with spherical objects $\A$.
We assume for simplicity that
\begin{myeq}[\label{eqgensphere}]
\A=\{\Sigma^{k}A_{0}\}_{k=0}^{\infty} 
\hsm\text{for some strict cogroup object}\hsm A_{0}~.
\end{myeq}

An \emph{elementary Stover object} in $\C$ is one of the form:
\begin{myeq}[\label{eqstov}]
B~:=~\colim~\left(A\xra{\inc} (C\Aj{j})_{j\in T}\right)~,
\end{myeq}
\noindent where \w[,]{A\in\A} and the colimit is of the diagram 
consisting of $A$, together with an inclusion \w{A\hra C\Aj{j}} into the cone 
on \w{\Aj{j}} (a copy of $A$) for each \w[.]{j\in T} The set $T$ is called the 
\emph{null set} for $B$. Note that $B$ is still in \w[,]{\CA} and is still a 
cogroup object in $\C$. 

A \emph{Stover object} is any coproduct \w{B=\coprod_{i\in I}~B\li{i}}
of elementary Stover objects \w[.]{\{B\li{i}\}_{i\in I}}

The \emph{Stover category}, denoted by \w[,]{\CAs} is the full
sub-simplicial category of \w{\sCA} consisting of all Stover objects
such that the cardinalities of the indexing set $I$ for the coproduct,
and of the null sets \w{T\li{i}} for each coproduct summand
\w[,]{B\li{i}} are bounded by a fixed limit cardinal $\kappa$ (see
Remark \ref{rcask} below).  

Evidently, \w{\CAs} is a small \ww{\sCA}-variant (\S \ref{dcap}). Any \Ama\ based 
on \w{\CAs} will be called a \emph{\Sma}, and the realizable \Sma\ for
any \w{Y\in\C} will be denoted by \w[.]{\fMAs Y} The category of all
\Sma s will be denoted by \w[.]{\MAs} 

Similarly, any \Ama\ based on the canonical minimal small \ww{\sCA}-variant \w{\CAm} 
(\S \ref{egcap}) will be called a \emph{minimal \Ama}, and the minimal \Ama\ for $Y$ 
will be denoted by \w[.]{\fMAm Y} 
\end{defn}

\begin{lemma}\label{lmapst}
For any \w[,]{Y\in\C} the mapping spaces of the \Sma\ \w{\fXs=\fMAs Y} are 
canonically determined by the minimal \Ama\ \w[.]{\fXm=\fMAm Y}
\end{lemma}

\bproof
For \w[,]{A\in\A} set \w[.]{\fXs\lin{A}:=\fXm\lin{A}}
If $B$ is an elementary Stover object as in \wref{eqstov} (with
\w[),]{T\neq\emptyset} we define \w{\fXs\lin{B}} to be the pullback in
\w[:]{\Sa} 
\mydiagram[\label{eqzeromaps}]{
\ar @{} [drr] |<<<{\framebox{\scriptsize{PB}}}
\fXs\lin{B} \ar[rr] \ar[d]_{f} && 
\prod_{j\in T} P\fXm\lin{\Aj{j}} \ar[d]^{\prod_{j}p_{j}}\\ 
\fXm\lin{A} \ar[rr]^<<<<<<<<{\Delta} && \prod_{j\in T} \fXm\lin{\Aj{j}} 
}
\noindent (where \w{P\fX} is the path functor of \S \ref{rloops} and
$\Delta$ is the diagonal). 

If \w{B=\coprod_{I\in I}~B\li{i}} is a coproduct of elementary Stover
objects, we set: 
\begin{myeq}[\label{eqcoprodst}]
\quad\hfill\fXs\lin{B}~:=~\prod_{I\in I}~\fXs\lin{B\li{i}}\hfill\Box
\end{myeq}

\begin{remark}\label{rmapst}
If $\fX$ is any fibrant \Ama\ based on the minimal small \ww{\sCA}-variant \w[,]{\CAm} 
we may use \wref{eqzeromaps} and \wref{eqcoprodst} to define the mapping spaces of the 
corresponding \Ama\ \w{\fXs} based on \w[.]{\CAs} Of course, this does not determine 
the action of \w{\CAs} on \w[.]{\fXs}
\end{remark}

\begin{lemma}\label{lrealistic}
If an \Ama\ $\fX$ based on \w{\CAm} is realistic (\S \ref{drealistic}), so is
corresponding \Ama\ \w{\fXs} based on \w[.]{\CAs} 
\end{lemma}

\begin{proof}
Since the right vertical map in \wref{eqzeromaps} is a fibration, so is 
\w[,]{f:\fXs\to\fXm} and \wref{eqzeromaps} is a homotopy pullback. 
Thus if $\zeta$ in \wref{eqloops} and $\theta$ in \wref{eqcoprodma} are weak 
equivalences for $\fX$ whenever the maps in \wref{eqrealistic} are, 
the same is true for \w[.]{\fXs} 
\end{proof}

\begin{cor}\label{crealistic}
Under assumption \wref[,]{eqgensphere} all the mapping spaces of a realistic 
\Sma\ \w{\fXs} are determined up to weak equivalence by the single 
simplicial set \w[.]{\fX\lin{A_{0}}}
\end{cor}

In the dual case (\S \ref{rdpa}), when we have homotopy group objects 
\w{\{W_{n}\}_{n=1}^{\infty}} in $\C$ with each \w[,]{W_{n}=\Omega W_{n+1}} 
it is not enough to know the single mapping space \w[;]{\map_{\C}(X,W_{1})} 
in this case we need its \ww{\Omega^{\infty}}-structure.

\begin{defn}\label{dstov}
Let $\fX$ be an \Ama\ based on a small \ww{\CA}-cvariant \w[,]{\hCA}
and \w[.]{B\in\hCA} For each \w{\phi\in \fX\lin{A}_{0}} we call the
pullback \w{N^{\phi}} (in \w[):]{\Sa} 
$$
\xymatrix@R=25pt{
\ar @{} [drr] |<<<<{\framebox{\scriptsize{PB}}}
N^{\phi} \ar[rr] \ar[d] && P\fX\lin{B} \ar[d]^{p}\\ 
\phi \ar[rr]^<<<<<<<<<<{\inc} && \fX\lin{B} 
}
$$
\noindent the \emph{space of nullhomotopies} for $\phi$. (It will be
empty if $\phi$ is not null-homotopic.)

If \w{\hCA} is any small \ww{\sCA}-variant containing $\A$ itself, and
$\fY$ is any \Ama\ based on \w[,]{\hCA} the \emph{Stover construction}
on $\fY$ is the Stover object given by: 
\begin{myeq}[\label{eqstovc}]
 K\fY~:=~\ \coprod _{A\in\A} \ 
     \coprod _{\phi\in \fY\lin{A}_{0}} \ \colim \left( A\li{\phi}\xra{\inc}
                  (CA\li{\Phi})_{\Phi\in N^{\phi}_{0}}\right)~.
\end{myeq}
\noindent This defines a functor \w[.]{K:\MAs\to\C} 
\end{defn}

%
%
\begin{prop}\label{pscm}
The composite \w{L:=K\circ\fMAs:\C\to\C} is a comonad on $\C$.
\end{prop}

\begin{proof}
Note that \w{K\fY} depends only on the $0$-simplices
  \w{\rz\fY:=(P\fY_{0}\to\fY_{0})} of the path fibration $\rho$ (\S
  \ref{rloops}). Because $\rho$ is a functor, any map of \Ama s
  \w{\Psi:\fY\to\fZ} induces a map of the indexing categories for the
  colimit \wref[.]{eqstovc} Again, this depends only on \w[.]{\rz\Psi}
  This in turn induces a map \w[.]{K\Psi:K\fY\to K\fZ} Thus we have
  defined a functor \w[.]{K_{0}:\rz\MAs\to\C} We show that the functor
  \w{K_{0}} is left adjoint to \w[:]{\rz \fMAs:\sC\to\rz \MAs} 

Given \w{f:K\fY\to X} in $\C$, we define \w{\hat{f}:\rz\fY\to\rz\fMAs X}
by sending \w{\phi\in\fY\lin{A}_{0}} to 
\w[,]{f\rest{A\li{\phi}}\in(\fMAs X)\lin{A}_{0}=\map_{\sC}(A,X)_{0}}
and similarly for \w{\Phi\in\fY\lin{A}_{1}} with \w{d_{0}\Phi=0} and 
\w{d_{1}\Phi=\phi} (using \wref[).]{eqconepath}

Conversely, given \w{\psi:\rz\fY\to\rz(\fMAs X)} in \w[,]{\rz\MAs}
we define \w{\tilde{\psi}:K\fY\to X} using the fact that \w{K\fY} is defined by
the colimit \wref[,]{eqstovc} so it is enough to define a map of diagrams,
given by \w{\psi(\phi):A\li{\phi}\to X} for \w{\phi\in\fY\lin{A}_{0}} and
\w{\psi(\Phi):CA\li{\Phi}\to X} for \w[,]{\Phi\in(P\fY\lin{A})_{0}}
again using \wref[.]{eqconepath}

Since we can factor \w{L:=K\circ\fMAs} as the composite
\w{K_{0}\circ(\rz\fMAs)} of an adjoint pair of functors, the functor
\w{L:\C\to\C} is a comonad (cf.\ \cite[\S 4]{BorcH2}).
\end{proof}

\begin{remark}\label{rrealistic}
Let \w{\fX} be a fibrant \Sma, and assume that each \w{A\in\A} is a strict 
cogroup object in $\C$.  Thus \w{\fX\lin{A}} is the underlying 
simplicial set of a simplicial group. Moreover, since the structure maps 
in \wref{eqzeromaps}and \wref{eqcoprodst} are all maps of simplicial groups 
(see \S \ref{rloops}), the same is true of \w{\fX\lin{B}} for \w[.]{B\in\CAs} 
(Of course, the composition maps in \w{\CAs} need not be hommomorphisms, so
$\fX$ is not necessarily enriched in $\G$.)

If \w{K\sp{e}} is the zero-component of \w[,]{K=\fX\lin{B}} we thus have two 
canonical short exact sequences of simplicial groups (resp., groups):
$$
\xymatrix@R=25pt{%
1 \to K\sp{e} \to K \to \pi_{0}K \to 1 &&
1 \to PK_{0} \ar[r] & K_{1}\sp{e} \ar[r]^<<<<{d_{0}} & 
K\sp{e}_{0} \to 1 \ar@/^1pc/[l]^{s_{0}} 
}
$$
\noindent This implies that \w{K_{1}} is canonically determined \emph{as a set} by
\w[,]{K_{0}} \w[,]{\pi_{0}K} and \w[.]{PK_{0}} In other words, \w{\rz\fX} and 
\w{\pi_{0}\fX} together determine \w{\csk{1}\fX} up to isomorphism 
(and of course conversely).
\end{remark}

\begin{defn}\label{dscm}
We call \w{L:\C\to\C} the \emph{Stover comonad} on $\C$. 

The counit \w{\var:L\to\Id} for $L$ is the ``tautological'' natural transformation 
\w[,]{\var_{X}:K(\fMAs X)\to X} which sends the copy of $A$ indexed by 
\w{\phi\in(\fMAs X)\lin{A}_{0}=\Hom_{C}(A,X)} in \wref{eqstovc} to $X$ by $\phi$, 
and similarly for the cones \w[.]{CA\li{\Phi}}

The comultiplication \w{\mu:L\to L^{2}} is induced by the natural inclusion 
\w[,]{\nu:K\fY\to K(\fMAs (K\fY))} defined for any \Sma\ $\fY$, 
which sends \w{A\li{\phi}} in \w{K\fY} identically to the copy of $A$ in 
\w{K(\fMAs (K\fY))} indexed by the inclusion \w[.]{A\li{\phi}\hra K\fY}

The \emph{Stover resolution} of an object \w{Y\in\C} is the simplicial resolution
\w{\Qd} of $Y$, where \w{Q_{n}:=L^{n+1}Y} for each \w{n\geq 0} (and the face and 
degeneracy maps are induced by $\eta$ and $\mu$).
\end{defn}

\begin{remark}\label{rcskk}
If we extend $K$ to a \emph{simplicial} functor \w[,]{\tK:\MAs\to\sC}
it factors through \w[,]{\bK:\rho\MAs\to\sC} so \w{\csk{n}\bK} depends
on \w[,]{\csk{n}\rho\MAs} which is determined in turn by \w[.]{\csk{n+1}\MAs}
\end{remark}

\begin{prop}\label{pstover}
If \w{s\C} is a strict \ww{E^{2}}-model category with spherical objects $\A$, 
the Stover resolution defines a one-to-one correspondence 
between objects \w{Y\in\C} up to $\A$-equivalence (\S \ref{dmodels}) and 
weak equivalence of simplicial objects \w{\Qd\in s\C} with \w{\piA\Qd\simeq\BL} 
(where \w[).]{\Lambda\cong\piA Y} 
\end{prop}

\begin{proof}
By \cite[\S 3.3]{DKStE} the simplicial object \w{\Qd} is cofibrant in
the resolution model category structure on \w[,]{s\C} and by 
\cite[\S 2]{StoV}, the map \w{\var:\Qd\to\co{Y}} (induced by $\eta$)
is a weak equivalence. Thus \w{\piA\Qd\simeq\BL} by \wref[.]{eqspiral}
From \wref{eqrealmap} we see that
\w[.]{\|\map_{\C}(A_{0},\Qd)\|\simeq\map_{\C}(A_{0},R(\Qd))} Applying
the Bousfield-Friedlander spectral sequence of \cite[Theorem B.5]{BFrieH} 
to the bisimplicial set \w[,]{\Md:=\map_{\C}(A_{0},\Qd)} with
\begin{myeq}[\label{eqbfried}]
E^{2}_{s,t}=\pi_{s}\pi_{t}\Md~\Rightarrow~\pi_{s+t}\|\Md\|~,
\end{myeq}
\noindent we conclude that \w[,]{\pis\|\Md\|\cong\Lambda} and thus
\w[.]{\piA(R\Qd)\cong\pis\map_{\C}(A_{0},R(\Qd))\cong\Lambda}
This is an isomorphism of \PAa s, since \w{\fMAs\Qd} is a simplicial
\ma, and so applying \w{\|-\|} to each bisimplical set \w{\fMAs\Qd\lin{A}}
\wb{A\in\A} yields a \ma, which is actually determined by
\w{\|\Md\|=(\|\fMAs\Qd\|)\lin{A_{0}}} by Corollary \ref{crealistic}.
Thus \w{R\Qd} is $\A$-equivalent to $Y$. 
Functoriality of the Stover construction (and of the spectral
sequence) shows that the correspondence of weak $\A$-homotopy types is
one-to-one. 
\end{proof}

\begin{remark}\label{rcask}
We can now explain how the cardinal $\kappa$ of \S \ref{dstovcat} is chosen\vsm: 

Given a \PAa\ $\Lambda$, the collection of all homotopy types of
objects \w{Y\in\ho\C} with \w{\Lambda\cong\piA Y} is a set (as can be
seen by considering all choices of $k$-invariants for cofibrant
replacements of \w{\co{Y}} in \w[).]{s\C} 

Define $\kappa$ to be the smallest limit cardinal such that each such 
homotopy type $Y$ has a Stover resolution in which each of the sets 
\w{\fMAs(L^{n}Y)\lin{A}_{0}} and \w{N^{\phi}} for
\w{\phi\in\fMAs(L^{n}Y)\lin{A}_{0}} in \wref[,]{eqstovc} for each
\w{A\in\A} and \w[,]{n\geq 0} has cardinality $\leq\kappa$. 
\end{remark}

\begin{mysubsection}{Extending the Stover comonad}\label{sescm}
Applying the functor \w{\fMAs} to the augmented simplicial object
\w{\Qd\to Y}  over $\C$ yields an augmented simplicial object
\w[.]{\fMAs\Qd\to\fMAs Y} We can think of this as coming from a
monad $\fL$ on \emph{realizable} \Sma s, given by
\w[,]{\fL(\fY):=\fMAs(\tK\fY)} with counit \w{\eta:=\fMAs(\var)} right
inverse to the unit \w{\xi:\fMAs Y\to \fMAs(K(\fMAs Y))} 
(sending \w{\phi:A\to Y} to the inclusion 
\w[).]{A\li{\phi}\hra K(\fMAs Y))} Because $\var$ was a counit for
$L$, the following square commutes: 
\mydiagram[\label{eqalgebra}]{
\fL\fL\fX \ar[r]^{\mu_{\fX}} \ar[d]_{\fL(\eta_{\fX})} & 
\fL\fX \ar[d]_{\eta_{\fX}} \\
\fL\fX \ar[r]_{\eta_{\fX}} & \fX
}
\noindent for \w{\fX=\fMAs Y} (cf.\ \cite[\S 4.1]{BorcH2}.

We observe that even though the simplicial functor \w{\fMA} does not usually
preserve coskeleta (even for \w{A=\bS{1}} in \w[),]{\Sa} 
we deduce from Remark \ref{rcskk} that:
\begin{myeq}[\label{eqcskl}]
\csk{n}\fL\fX \hsm\text{is determined by}\hsm \csk{n+1}\fX
\end{myeq}
\noindent because \w{\tK} actually lands in \w[,]{\CAs} so $\fL$ takes
values in free \Sma s (\S \ref{dfreema}).
\end{mysubsection}

\begin{defn}\label{dtalg}
A fibrant \Sma\ $\fX$ is called an $\fL$-\emph{algebra} if it is equipped 
with a splitting \w{\eta_{\fX}:\fL\fX\to\fX} for 
\w[,]{\xi:\fX\to\fL\fX} such that \wref{eqalgebra} commutes.
\end{defn}

%
%
\begin{prop}\label{ptalg}
Any realistic \Sma\ $\fY$ can be realized, up to $\A$-equivalence.
\end{prop}

\begin{proof}
Iterating the functor $\fL$ on $\fY$ yields an augmented simplicial \Sma\ 
\w[,]{\fVd\to\fY} and since \w[,]{\fL=\fMAs\circ K} in fact 
\w[.]{\fVd=\fMAs\Qd} Here \w{\Qd} is the simplicial Stover object with 
\w{Q_{0}=K\fY} and \w{Q_{n}:=K\fV_{n-1}} for \w[.]{n\geq 1} The extra 
face map \w{d_{n}:Q_{n}\to Q_{n-1}} is 
\w[.]{K\fL^{n-1}(\eta):K\fL^{n-1}\fMAs\K\fY\to K\fL^{n-1}\fY}
where \w{\eta:\fMAs\K\fY\to \fY} is the $\fL$-algebra structure map.

Proposition 2.6 of \cite{StoV} shows that if \w{\Qd\in s\C} is 
the Stover resolution of \w[,]{X\in\C} then \w{\pis\Qd} is a free
\PAa\ resolution of \w[.]{\piA Y} The proof does not in fact depend on
the existence of $X$, but only on its \ma\ \w[.]{\fY:=\fMAm X} 
Here we use the fact that $\fY$ is realistic. 
Thus we deduce that the simplicial \PAa
\w{\Gd:=\piA\Qd\cong\pi_{0}\fVd} is a free \PAa\ resolution of
\w[.]{\Lambda=\pi_{0}\fY} Thus the spectral sequence of
\wref{eqbfried} collapses, showing that \w{R\Qd} realizes $\Lambda$. 

Finally, by combining the weak equivalences of \S \ref{drealistic}
with Lemma \ref{lmapst}, we deduce that \w{\|\fVd\|} (the realization
functor applied to each simplicial space \w{\fVd\lin{B}} for
\w[)]{B\in\CAs} is an \Sma, which is weakly equivalent to
\w[,]{\fMAs(R\Qd)}as well as to the original \ma\ $\fY$.
\end{proof}

Since every realizable \Sma\ is realistic, this shows:

\begin{cor}\label{ctalg}
The correspondence of Proposition \ref{pstover} actually factors through 
the category of realistic $\fL$-algebras, up to weak equivalence.
\end{cor}

%
%
\sect{Realizing mapping algebras}
\label{crma}

In order to solve the first problem mentioned in the beginning of 
Section \ref{cma}, we must reinterpret the inductive approach to 
realizing a \PAa\ $\Lambda$ described in Section \ref{crpa} as an 
inductive process for realizing \ma s. For this we need:

\begin{defn}\label{dnralg}
A map of \Ama s \w{\fff:\fX\to\fY} is called an $n$-\emph{equivalence} if 
\w{\Po{n}\fff} is a weak equivalence of \Ama s. Similarly, a map \w{f:X\to Y}
in $\C$ is called an $n$-$\A$-\emph{equivalence} if \w{\fMA f} is an
$n$-equivalence of \Ama s.

A $\fL$-algebra $\fX$ is called a \emph{$n$-realistic $\fL$-algebra} if:
 
\begin{enumerate}
\renewcommand{\labelenumi}{(\alph{enumi})}
\item \w{r\q{n}:\fX\to\Po{n}\fX} is a weak equivalence of \Ama s.
\item The map $\zeta$ in \wref{eqloops} is an \ww{(n-1)}-equivalence in \w{\Sa}
      whenever the first map in \wref{eqrealistic} is a weak equivalence.
\item The map $\theta$ in \wref{eqcoprodma} is an $n$-equivalence 
      whenever the second map in \wref{eqrealistic} is a weak equivalence.
\end{enumerate}
\end{defn}

\begin{remark}\label{rnralg}
Note that we cannot expect to do better than (b) above, since
\w{\Omega P^{n}K} is just \w{P^{n-1}\Omega K} for any \w[.]{K\in\Sa}
Thus even under Assumption \wref[,]{eqgensphere} where for a realistic
\Sma\ $\fX$, the simplicial set \w{\fX\lin{A_{0}}} determines
\w{\fX\lin{B}} for any $B$ in \w{\CAs} up to weak equivalence, in the
$n$-realistic case \w{\fX\lin{\Sigma A_{0}}} carries more information
than \w{\Omega\fX\lin{A_{0}}} does. 
\end{remark}

We can now refine Corollary \ref{ctalg} as follows:

%
%
\begin{prop}\label{pstrefined}
There is a one-to-one correspondence between $n$-realistic
$\fL$-algeb\-ras $\fX$ with \w{\pi_{0}\fX\cong\Lambda\in\PAAlg} and
$n$-th quasi-Postnikov sections for $\Lambda$, up to weak equivalence.
\end{prop}

\begin{proof}
Let $\fX$ be an $n$-realistic $\fL$-algebra, so its structure map 
\w{\eta=\eta_{\fX}} factors through 
\w[.]{\Po{n}\fL\fX\to\fX\simeq\Po{n}\fX=\csk{n+1}\fX} 
We wish to construct the Stover resolution \w{\fVd\to\fX} as in the proof
of Proposition \ref{ptalg}. For all \w[,]{k\geq 0} the objects 
\w{\fV_{k}:=\fL^{k+1}\fX} depend only on \w[,]{\rz\fX} which is determined by 
\w[.]{\Po{0}\fX=\csk{1}\fX} Similarly, all the degeneracy and face maps, in
all simplicial dimensions,  are determined by \w[,]{\fV_{0}\in\CAs} except 
for \w[,]{d_{k}:\fV_{k}\to\fV_{k-1}} which is \w[.]{\fL^{k}\eta} 
By \wref[,]{eqcskl} this map itself, as an arrow in in \w[,]{\uCAs\subseteq\C} 
depends only on \w[.]{\csk{k}\eta:\csk{k}\fV_{0}\to\csk{k}\fX} 
Thus $\eta$ determines the \ww{n+1}-st truncation \w{\tau_{n+1}\fVd} of 
\w[,]{\fVd} and thus \w[.]{\Po{n}\fVd}

Conversely, if we can construct \w{\tau_{n+1}\fVd} for $\fX$, this is
equivalent (as in the proof of Proposition \ref{ptalg}) to
constructing \w{\tau_{n+1}\Qd} for the (putative) object \w{Y\in\C}
realizing $\fX$, with \w[.]{\tau_{n+1}\fVd:=\fMAs(\tau_{n+1}\Qd)}
Thus we have an $n$-th quasi-Postnikov section for $\Lambda$, which we
denote by \w{\Qnd{n}\in s\C} (see Definition \ref{dquap}). Applying
the $n$-th Postnikov section functor \w{\Po{n}:\tau_{n+1}s\Sa\to s\Sa}
to each \ww{(n+1)}-truncated simplicial set \w{\tau_{n+1}\fVd\lin{A}}
yields the corresponding quasi-Postnikov section \w[,]{\fVnd{n}\in
  s\MAs} with \w[.]{\fVnd{n}:=\fMAs(\Qnd{n})} This is because each
\w{\Qn{n}{k}} for \w{k\geq n+2} is constructed as a matching object
(cf.\ \cite[\S 2.1]{DKStB}), which is a limit, so it commutes with
mapping spaces. 

In particular, \w{\piA\Qn{n}{k}\cong(\pi_{0}\fVn{n}{k}\lin{A})_{A\in\A}} 
for all \w[.]{k\geq 0} Thus from \wref{eqethg} we see:
\begin{myeq}[\label{eqqpo}]
\pi_{k}\pi_{0}\fVnd{n}\lin{A}~=~\pi_{k}(\piA\Qnd{n})\lin{A}\cong\begin{cases}
		\Lambda\lin{A} &~~\text{for \ } k=0\\
		(\Omega^{k+1}\Lambda)\lin{A} &~~\text{for \ } k=n+2,\\
		 0     &~~\text{otherwise}~\end{cases}
\end{myeq}
\noindent for any \w{A\in\A} (and thus, using Lemma \ref{lmapst}, for
any \w[).]{B\in\CAs} The Bousfield-Friedlander spectral sequence
\wref{eqbfried} for the bisimplicial set \w{\fVnd{n}\lin{A}} converges
to \w[,]{\pis\|\fVnd{n}\|} because \w{\pi_{0}\fVnd{n}\lin{A}} is a
simplicial group. Moreover, the first possible differential is
\w[,]{d^{n+2}:(\Omega^{n+1}\Lambda)\lin{A})_{0}\to\Lambda\lin{A}_{n+1}} so
\w{\pi_{i}\|\fVnd{n}\lin{A}\|\cong\Lambda\lin{A}_{i}} for \w[.]{i\leq n} 
By naturality we deduce that that the map of \Sma s \w{\|\fVnd{n}\|\to\fX} 
is an $n$-equivalence, so \w[.]{\Po{n}\|\fVnd{n}\|\simeq\fX} 

In summary, each of \w[,]{\fVnd{n}\in s\MAs} \w[,]{\Qnd{n}\in s\C} and the 
$n$-realistic \Sma\ $\fX$ determines the other two.
\end{proof}

\begin{remark}\label{rstrefined}
Note that from the quasi-Postnikov section \w{\Qnd{n}\in s\C} we can
also recover an object \w{\Zn{n}:=R(\Qnd{n})\in\C} (using Definition
\ref{demc}), and we see that 
\w{[\Sigma^{i}A_{0},\Zn{n}]_{\C}\cong\pi_{i}\fX\lin{A_{0}}\cong
\Lambda\lin{\Sigma^{i}A_{0}}}
for \w[,]{0\leq i\leq n} by \wref[,]{eqrealmap} since \w{\Qnd{n}} is
$\A$-cofibrant and $\A$ is generated by \w{A_{0}} by \wref[.]{eqgensphere}

However, we can do more than this, by Remark \ref{rnralg}: the
inclusion of the subcollection of spherical objects 
\w{\An{k}:=\{\Sigma^{k}A_{0},\Sigma^{k+1}A_{0},\dotsc\}} in $\A$ induces a 
forgetful functor \w{\MA\to\MAn{k}} (which omits the simplicial set
\w{\fX\lin{\A_{0}^{i}}} \wb{0\leq i<k} from \w[).]{\fX} If we denote this
by \w[,]{\fX\mapsto\fX\q{k}} applying the procedure described in the proof of
Proposition \ref{pstrefined} to the \Ama\ \w{\fX\q{k}} (which is still
$n$-realistic) yields a new simplicial object \w[,]{\Qnd{n}\q{k}\in s\C} 
and \w{\Zn{n}\q{k}:=R\Qnd{n}\q{k}} now realizes the \PAa\
\w{\Omega^{k}\Lambda} through degree $n$. Moreover, there is a natural
\ww{(n-1)}-$\A$-equivalence \w{\Omega \Zn{n}\q{k+1}\to\Zn{n}\q{k}} for
each \w[,]{k\geq 0} induced by the maps $\zeta$ of \wref[.]{eqloops}

The collection of objects \w[,]{\{\Zn{n}\q{k}\}_{k=0}^{\infty}}
equipped with these structure maps, thus form an $n$-\emph{stem}, in
the sense of \cite{BBlaS}. In the case when \w{\C=\Ta} and
\w[,]{\A=\{\bS{k}\}_{k=1}^{\infty}} these behave like the collection
\w{\{\Po{n+k}X\lra{k-1}\}_{k=1}^{\infty}} of \ww{(k-1)}-connected  
covers of \ww{(n+k)}-Postnikov sections of a (putative) space $X$.
\end{remark}

\begin{lemma}\label{lcorhg}
Let $\fX$ be an $n$-realistic \Ama, and let \w{\Qnd{n}} be the
$n$-quasi Postnikov section in \w{s\C} corresponding to $\fX$ under
Proposition \ref{pstrefined}, with \w[.]{\\Lambda:=\pi_{0}\fX} Then
there is a natural isomorphism \w{\pi_{i}\fX\cong\pinat{i}{\Qnd{n}}}
as $\Lambda$-modules for all \w[.]{i\geq 0} 
\end{lemma}

\begin{proof}
From \wref[,]{eqqpo} \wref[,]{eqpnpsys} and \wref[,]{eqbfried} we see that
$$
\pi_{k}\fX\cong\pinat{k}{\Qnd{n}}\cong
\begin{cases}\Omega^{k}\Lambda &~~\text{for \ } k\leq n\\
		 0     &~~\text{otherwise~.}~\end{cases}
$$
\noindent To describe the natural identification, note that by 
Proposition \ref{pstrefined} we know that $\fX$ is $n$-equivalent to 
\w[.]{\|\fVnd{n}\|=\|\fMAs(\Qnd{n})\|} 
Since we assumed each \w{A\in\A} was a strict cogroup object in $\C$,
\w{K:=\|\fMAs(\Qnd{n})\|\lin{A}} has the natural structure of a
simplicial group. Therefore, an element in \w{\pi_{k}\fX\lin{A}} may
be represented by a Moore $k$-cycle $\phi$ in 
$$
Z_{k}K\subseteq K_{k}=\map_{\C}(A,\Qn{n}{k})_{k}=
\Hom_{\C}(A\otimes\Del[k],\Qn{n}{k})\subseteq
\Hom_{s\C}(\co{A}\hotimes\Del[k],\Qnd{n})
$$
\noindent (see \S \ref{snac}).

On the other hand, by \cite[Proposition 5.8]{DKStB} we can represent
an element of \w{\pinat{k}{\Qnd{n}}\lin{A}} by an element in
\w{\Hom_{\C}(A,Z_{k}\Qnd{n})} \wh that is by a map\w{f:A\to \Qn{n}{k}}
such that \w{d_{i}f=\ast} for all \w[.]{0\leq i\leq k} If
\w{\delta_{k}\in\Del[k]_{k}} is the non-degenerate $k$-simplex of
\w[,]{\Del[k]} we define
\w{\hat{\phi}:A\hotimes\Del[k]\to\Qnd{n}} by sending 
\w{A\otimes\{\delta_{k}\}} to \w{\Qn{n}{k}} by $f$, and extend by zero
to the other non-degenerate simplices of \w[.]{\Del[k]}
\end{proof}

\begin{defn}\label{daspalg}
For any \Ama\ $\fX$, the \emph{associated simplicial \PAa}
\w{\fPAd{\fX}} is defined by requiring \w{\fPA{\fX}{n}} to be the
\PAa\ induced by the action of \w{\pi_{0}\hCA} on each set of
$n$-simplices \w{\fX\lin{A}_{n}} of \w[.]{\fX\lin{A}\in\Sa} Note that
\w{\fPAd{\fX}} is itself an \Ama, and the quotient map
\w{h:\fX\to\fPAd{\fX}} is a map of \Ama s. 
\end{defn}

For simplicity, let us denote the cofibrant object \w{\Qd\in s\C} associated 
by Proposition \ref{ptalg} to a realistic \Ama\ $\fX$ by \w[.]{\Qnd{\infty}}

\begin{lemma}\label{lhurewicz}
Assume that $\fX$ is an $n$-realistic \Ama, with \w{0\leq n\leq\infty}
and \w{\Qnd{n}} is the object associated to $\fX$ by Proposition \ref{ptalg} 
(respectively, Proposition \ref{pstrefined}). There is a natural
isomorphism of \w{\fPAd{\fX}} with the simplicial \PAa\
\w[,]{\piA\Qnd{n}} and \w{h:\fX\to\fPAd{\fX}} induces the Hurewicz
homorphism \w{h_{\#}:\pinat{\ast}{\Qnd{n}}\to\pis\piA\Qnd{n}} of
\wref[.]{eqspiral} 
\end{lemma}

\begin{proof}
As in the proof of Lemma \ref{lcorhg}, we may replace $\fX$ by the $n$-equivalent 
\Ama\ \w[,]{\|\fMAs(\Qnd{n})\|} so that any element in\w{\fX\lin{A}_{k}} may be 
identified with a map \w[.]{f:A\otimes\Del[k]\to \Qn{n}{k}} 

Since \w[,]{Q:=\Qn{n}{k}\in\CAs} we may identify this with 
\w[,]{f^{\ast}(\Id_{Q})} for \w[,]{f\in(\fMAs Q)\lin{A}} so by definition of $h$
we have \w[.]{h(f)=[f]^{\ast}h(\Id_{Q})=[f]\in[A,Q]_{\C}=(\piA\Qn{n}{k})\lin{A}}
This identifies \w{\fPAd{\fX}} with \w[.]{\piA\Qnd{n}} From the description of
the Hurewicz homomorphism in \cite[\S 5]{DKStB}, we see that it coincides 
with \w[.]{h_{\#}}
\end{proof}

%
%
\begin{prop}\label{pobstrnr}
If \w{\fX=\fXn{n}} is an $n$-realistic \Ama\ with 
\w[,]{\Lambda:=\pi_{0}\fX\in\PAAlg} the obstruction to extending $\fX$ to
an \ww{(n+1)}-realistic \Ama\ \w{\fXn{n+1}} (with \w[)]{\Po{n}\fXn{n+1}=\fXn{n}}
is the \ww{(n+1)}-st $k$-invariant for \w[,]{\fPAd{\fX}} i.e.,
\w[.]{\tk_{n+1}\in\HAQ{n+3}{\Lambda}{\Omega^{n+1}\Lambda}}
\end{prop}

\begin{proof}
Again, let \w{\Qnd{n}} be the $n$-th quasi-Postnikov section for $\Lambda$
corresponding to $\fX$ under Proposition \ref{pstrefined}, 
with \w[.]{\fX\simeq\Po{n}\|\fMAs(\Qnd{n})\|}
By Lemma \ref{lhurewicz} and \wref{eqethg} we know \w{\fPAd{\fX}} has only
two non-zero homotopy groups: \w{\pi_{0}\fPAd{\fX}\cong\Lambda} and
\w[.]{\pi_{n+2}\fPAd{\fX}\cong\Omega^{n+1}\Lambda}

If the extension \w{\fXn{n+1}} exists, the fibration 
\w{p\q{n+1}:\fXn{n+1}\to\Po{n}\fXn{n+1}\simeq\fXn{n}} induces
\w[,]{p\q{n+1}_{\#}:\fPAd{\fXn{n+1}}\to\fPAd{\fXn{n}}} which is the identity on
\w{\pi_{0}\fXn{n+1}=\pi_{0}\fXn{n}\cong\Lambda} (again by Lemma \ref{lhurewicz}).

Since \w[,]{\fPAd{\fXn{n}}\simeq\Po{n+1}\fPAd{\fXn{n}}}
\w{p\q{n+1}_{\#}} factors via \w[,]{\Po{n+1}\fPAd{\fXn{n+1}}=\BL} so
the structure map  \w{\tp\q{n+2}:\fPAd{\fX}\to\Po{n+1}\fPAd{\fX}= \BL}
has a section $s$. This is equivalent to the vanishing of the
\ww{(n+1)}-st $k$-invariant
\w{\tk_{n+1}\in\HAQ{n+3}{\Lambda}{\Omega^{n+1}\Lambda}} for
\w[.]{\fPAd{\fX}}

Conversely, if the $k$-invariant \w{\tk_{n+1}} for
\w{\fPAd{\fX}\cong\piA\Qnd{n}} vanishes, then \w{\Qnd{n}} extends to
an \ww{(n+1)}-st quasi-Postnikov section \w{\Qnd{n+1}} for $\Lambda$,
by Theorem \ref{treal}, so we obtain
\w[,]{\fXn{n+1}:=\|\fMAs(\Qnd{n+1})\|} as required by Proposition 
\ref{pstrefined}. 
\end{proof}

\begin{remark}\label{rhurewicz}
Since the quotient map $h$ of \S \ref{daspalg} is surjective, and 
\w{\fPAd{\fX}\lin{B}} has the underlying structure of a simplicial
group for each \w[,]{B\in\CAs} $h$ is a fibration in
\w[.]{\MAs\subseteq\SaOC} In fact, we may identify the long exact
sequence in \w{\pis} for the fibration $h$ with the spiral exact
sequence \wref[,]{eqspiral} up to a re-indexing. 

If we denote the fiber of $h$ by \w[,]{B'\fX} we deduce from
\wref{eqpnpsys} and \wref{eqethg} that:
$$
\pi_{i}(\B'\fX)\cong\begin{cases}
\Omega^{i}\Lambda &~~\text{for \ }1\leq i\leq n+1\\
		 0     &~~\text{otherwise}~\end{cases}
$$
\noindent Looping back the fibration sequence for $h$, for each
\w{A\in\A} we obtain 
$$
\fX\lin{\Sigma A}=\Po{n}\fX\lin{\Sigma A} \xra{p\q{n}}
\Po{n-1}\fX\lin{\Sigma A}\xra{\zeta'} \Omega\fX\lin{A}=
\Po{n-1}\Omega\fX\lin{A}\xra{k'} E(\Omega^{n+1}\Lambda,n+1),
$$
\noindent where \w{\zeta'} is the weak equivalence of \S \ref{dnralg}(b), and 
\w{k'} is the (looped) \ww{(n-1)}-th $k$-invariant for \w[.]{\fX\lin{\Sigma A}}

We can think of \w{\zeta'} as the structure map for the $n$-stem
\w[,]{\Po{n}\fX} which is classified by \w[.]{h:\fX\to\fPAd{\fX}} If
we could produce a map \w{q:\fPAd{\fX}\to E(\Omega^{n+1}\Lambda,n+2)}
which is a \ww{\pi_{n+2}}-isomrohism, then \w{q\circ h} would be the
$n$-th $k$-invariant for \w[,]{\fX=\fXn{n}} which would define an
\ww{(n+1)}-realistic \Ama\ \w[,]{\fXn{n+1}} and thus an
\ww{(n+1)}-quasi-Postnikov section for $\Lambda$. 

Now the inclusion of the homotopy fiber of $q$ is a map 
\w{s:\BL\to\fPAd{\fX}} which is a section for 
\w[.]{p=p\q{n+2}:\fPAd{\fX}\to\Po{n+1}\fPAd{\fX}= \BL} Moreover, $q$ is then 
\w{\Po{n+2}} applied to the pinch map to of the cofiber of $s$, so that the 
existence of $q$ is equivalent to the existence of a section $s$ for $p$.
Both are equivalent as above to the vanishing of the the \ww{(n+1)}-st
$k$-invariant for \w[.]{\fPAd{\fX}} 

Thus we can interpret this $k$-invariant, in the context of stems, as the 
obstruction to gluing the $n$-windows of an $n$-stem to produce an 
\ww{(n+1)}-stem. 
\end{remark}

We can thus summarize the results of this section in the following:

%
%
\begin{thm}\label{tobsma}
Let $\C$ be a strict \ww{E^{2}}-model category, with spherical objects $\A$
satisfying \wref[,]{eqgensphere} and let $\Lambda$ be a \PAa.
Proposition \ref{pobstrnr} then provides an inductively defined sequence of 
\ww{\SaO}-cohomology classes \w{\tk_{n}\in\HSO{n+2}{\BL}{\Omega^{n}\Lambda}}
\wb{n=1,2,\dotsc} to producing a realistic \Sma\ realizing $\Lambda$ \wh
which is equivalent to realizing $\Lambda$ in $\C$.\hfill$\Box$
\end{thm}

\begin{remark}\label{robsma}
We now interpret the classes \w{\tk_{n}} in the context of the Dwyer-Kan-Smith 
theory of Section \ref{cdoc}. The homotopy-commutative diagram which 
we are trying to rectify will be indexed by the $\Lambda$-pair 
\w[,]{(\D,\Dpp):=(\CAstp,\CAst)} defined as in Example \ref{egpostma}.
As in \S \ref{srelvo}, identifying the (ordinary) category \w{\CAst} as the 
zero-simplices of the (simplicially enriched) \w{\CAs} and applying degeneracies
gives the required simplicial map \w[.]{\tX_{0}:\co{\CAst}\to\CAs} 

However, the obstruction theory of Section \ref{cdoc} does not quite 
apply in our situation, since to begin with we do not have given a
simplicially enriched category \w{\CX} (\S \ref{rcx}) extending
\w{\CAs} \wh \emph{scil.} a (hopefully realistic) \Sma\
$\fX$. Instead, we construct $\fX$ by induction on its quasi-Postnikov
tower \w{(\fXn{n})_{n=0}^{\infty}} of $n$-realistic \Sma s. At the
beginning of the process we can always choose a $0$-realistic \Sma\
\w{\fXn{0}} realizing $\Lambda$, as well as an extension 
\w{X_{0}:\Po{0}\CAstp\to\fXn{0}} for the given \w[.]{\tX_{0}}

In view of Proposition \ref{pobstrnr}, we do not actually need to lift
\w{X_{0}} to the successive $n$-realistic \Sma s \w[,]{\fXn{n}}
but only to their \PAa\ versions \w{\fPAd{\fXn{n}}} These are \Sma s, though they 
are not $n$-realistic. Moreover, the Dwyer-Kan-Smith obstructions 
of \S \ref{srelvo} reduce in this case to the $k$-invariants \w[,]{\tk_{n}} 
as in the proof Proposition \ref{pobstrnr}.
\end{remark}

\end{document}